\documentclass{article}
\usepackage{amssymb,amsmath,amsthm}
\usepackage{bm}

\usepackage{graphicx}

\def\bH{{\bm H}}
\def\bQ{{\mathbb Q}}
\def\bR{{\mathbb R}}
\def\bS{{\mathbb S}}
\def\bT{{\mathbb T}}
\def\bZ{{\mathbb Z}}
\def\cA{{\cal A}}
\def\cC{{\cal C}}
\def\cD{{\cal D}}
\def\cM{{\cal M}}
\def\cN{{\cal N}}
\def\cP{{\cal P}}
\def\cS{{\cal S}}
\def\cZ{{\cal Z}}
\def\mf{{\mathfrak f}}
\def\mg{{\mathfrak g}}
\def\mh{{\mathfrak h}}
\def\Rt{\bR^3}

\def\Tt{\bT^3}
\def\Zt{\bZ^3}
\def\irr{\hbox{\rm irr}}
\def\span{\hbox{\rm span}}
\def\soul{\hbox{\rm soul}}
\def\Star{\hbox{\rm Star}}
\def\RPd{\bR\hbox{\rm P}^2}
\def\RPt{\bR\hbox{\rm P}^2}
\def\oN{\overline{\cN}}
\def\oM{\cM}

\def\rmU{\uppercase\expandafter{\romannumeral1}}
\def\rmD{\uppercase\expandafter{\romannumeral2}}
\def\rmT{\uppercase\expandafter{\romannumeral3}}

\def\fig#1#2{
  \begin{figure}
    \begin{center}
     \includegraphics[width=15cm,bb=140 0 800 734]{#1}
    \end{center}
   \caption{#2}
   \label{fig:#1}
  \end{figure}
}

\newtheorem{definition}{Definition}
\newtheorem{theorem}{Theorem}

\newtheorem{corollary}{Corollary}
\newtheorem{lemma}{Lemma}

\title{Topology of plane sections of periodic polyhedra with an application to the Truncated Octahedron}

\author{Roberto De Leo\thanks{INFN $<$roberto.deleo@ca.infn.it$>$, Dept. Physics, U. of Cagliari and Dept. of Mathematics, U. of Cagliari $<$deleo@unica.it$>$, Cagliari, Italy.}}

\begin{document}
\maketitle

\begin{abstract}
The main results of A. Zorich and I. Dynnikov about plane sections of periodic surfaces
are extended to the PL case. As an application, the Stereographic Map of a truncated
octahedron, extended to the whole $\Rt$ by periodicity, is analyzed numerically.
\end{abstract}



\pagestyle{myheadings}
\thispagestyle{plain}
\markboth{Roberto De Leo}{Plane sections of a periodic polyhedron}

\section{Introduction}

The problem of the asymptotics of plane sections of smooth periodic surfaces, extracted
from Physics literature by S.P. Novikov in 1982~\cite{Nov82}, turned out to be much richer
then expected, leading ultimately to the association of a fractal on $\RPt$ to every
element of a large class of triply periodic functions in $\Rt$ (see Section 2).
\par
In order to visualize for the first time Stereographic Maps associated to triply periodic
smooth functions and to get numerical confirmations of a Novikov conjecture, claiming that
the Hausdorff dimension of such fractals is stricly between 1 and 2, we developed a C++
library and used it to investigate two smooth cases \cite{DL03a}. Unfortunately, the
running time of our numerical explorations grows way too much as soon as we sample with
resolutions big enough to get a hint of the fractals, mainly because the number of
polygons of the meshes approximating curved smooth surfaces (e.g. to retrieve its
intersection with a plane) gets soon very big.
\par
This fact suggests that, from the numerical point of view, polyhedra are the best surfaces
to study, at the very least for the obvious reason that the number of polygons needed to
describe them at any resolution is constant. Moreover, such constant can be rather small
even in non trivial cases, e.g. like in the case of the ``extended truncated octahedron'',
presented in this paper, which has just eight hexagonal faces.
\par
Since the main results of the theory, due to A.V. Zorich~\cite{Zor84} and
I.A. Dynnikov~\cite{Dyn97,Dyn99}, refer to the case of smooth surfaces only, in Section 2
we will provide independent proofs for those theorems that make use of the smooth
structure and will mention the main properties of the system.  Then, in Section 3, we
present the algorithm we implemented to explore numerically the problem and the results
obtained in case of a the polyhedron obtained extending by periodicity the truncated
octahedron.
\section{Fundamental objects and theorems}
\subsection{Critical Points of height functions in polyhedra}
An analog of the Morse theory for height functions on polyhedra has been introduced by
T. Banchoff in \cite{Ban67,Ban70}. Here we recall the concepts relevant for the present
paper, slightly modified to cover the case of periodic polyhedra.
\begin{definition}
By ``embedded polyhedron'' $M\subset\Rt$ we mean a countable collection of cells
$K=\{C^r\subset\Rt\}_{r=0,1,2}$, where 0-cells are points (vertices), 1-cells are closed
connected segments and 2-cells are convex closed plane polygons, such that:
\par
\begin{enumerate}
    \item[{\textbf P1}] the boundary of any cell is union of cells of lesser degree;
    \item[{\textbf P2}] every cell having points in common with a higher degree cell is completely 
	contained inside it;
    \item[{\textbf P3}] K is locally finite, i.e. every vertex has a neighborhood in $\Rt$ that intersects 
	only finitely many cells;
    \item[{\textbf P4}] for any point $p\in M$, the union $\Star(p)$ of all cells containing that point
	is omeomorphic to an open disc.
\end{enumerate}
A triply periodic polyhedron is a polyhedron that is invariant with respect to a rank-3 
discrete subgroup $\Gamma\simeq\Zt$ of $\Rt$. 
Finally, by polyhedron $\oM\subset\Tt$ we mean the 
quotient $M/\Gamma\subset\Rt/\Gamma\simeq\Tt$ of a triply periodic polyhedron.
\end{definition}
\par
Every triply periodic polyhedron $M$ embedded in $\Rt$ is the lift of a compact polyhedron
$\oM$ embedded in $\Tt$. Since we are going to study polyhedra's plane foliations, we are
interested in the reciprocal relation between polyhedra and height functions (or,
equivalently, constant 1-forms):
\begin{definition}
A height function $h(p)=h^\alpha p_\alpha$ is called ``generic'' for the polyhedron $M$
if no edge of $M$ is perpendicular to the direction $\bH=(h^\alpha)$. Equivalently, a constant
1-form $\omega=h^\alpha dp_\alpha$ in $\Rt$ (resp. $\Tt$) is generic for $M$ (resp. $\oM$)
if no edge of $M$ (resp. $\oM$) is contained in a single leaf %
\footnote{Frobenius theorem grants that the distribution $\omega=0$ is integrable iff the 1-form
$\omega$ is closed.  In this case, the leaves induced by $\omega=h^\alpha dp_\alpha$ in
$\Tt$ are the projections of the $\Rt$ planes perpendicular to $\bH=(h^\alpha)$; those induced on a
polyhedron $\oM$ are the intersections of these leaves with the polyhedron.}%
of $\omega$.
\end{definition}
\par
Since height functions are not single-valued on $\Tt$ (unless $\bH$ is an integer
direction, i.e. parallel to a lattice vector) while their differentials $\omega=dh$ are
always well-defined in both $\Rt$ and $\Tt$, we will refer mostly to 1-forms from now on.
From the definition above it is clear that, like in the smooth case, the set of
non-generic 1-forms has zero measure.
\par
Note that the foliation induced on $\Tt$ (and therefore on $\oM$) by $\omega$ does not
change by multiplying the 1-form by a non-zero scalar, so from now on we will think of
$\omega$ at the same time as a constant 1-form and as a point in $\RPt$.
%
\begin{definition}
The index of a point $p\in M$ (or, equivalently, $[p]\in\oM$) with respect to a
generic constant 1-form $\omega$ is the integer $i(p,\omega)=1-s/2$, where $s$ is the
number of segments, having $p$ as one of their extremes, in which the leaf of $\omega$
passing through $p$ cuts $\Star(p)$.
If $i(p,\omega)=0$, i.e. if the $\Star(p)$ is cut in exactly two components, the point is
said ``regular''; otherwise it is called ``critical''.
\end{definition}
\par
\begin{figure}
\begin{center}
\includegraphics[width=5cm]{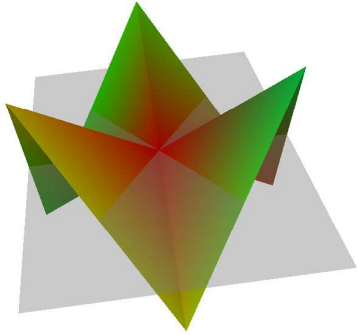}
\end{center}
\caption{%
Piece-wise linear Monkey Saddle: in the smooth case, a however small perturbation would be
enough to resolve this critical point in a pair of elementary saddles; in the piece-wise
linear case instead it is stable and, therefore, generic.
}
\label{fig:ms}
\end{figure}

The set of critical points for any generic constant 1-form is of course a subset of the
set of vertices; exactly as in the smooth case, minima and maxima have index $+1$ and
non-degenerate saddles have index $-1$.
The main difference between the smooth and the piece-wise linear (PL) case, for our
purposes, is that saddles that are unstable in the former case, i.e. disappear for small
perturbations of the 1-form direction, are stable in the latter and therefore cannot be
disregarded; the simplest example is provided by the ``Monkey saddle'', which has index
$-2$ (see Figure~\ref{fig:ms}).
\par
Nevertheless, an analog of the critical point theorem for generic 1-forms holds:
\begin{theorem}
If $M$ is a triply periodic polyhedron, invariant by the action of the rank-3 group
$\Gamma\subset\Zt$, and $\omega$ is generic for $M$, then
$\sum_{[p]\in\oM}i(p,\omega)=\chi(\oM)$, where $[p]$ is the set of vertices
$\Gamma$-equivalent to $p$ and the sum is hence extended to any set of inequivalent
vertices. Equivalently, if we set to 1 the volume of a Dirichlet domain of $\Gamma$,
then the average of the Euler characteristic of M converges to the Euler characteristic of
$\oM$:
$$\overline{\chi(M)}=\displaystyle\lim_{R\to\infty}\sum_{\|p\|<R}i(p,\omega)/Vol(B_R)=\chi(\oM)$$
\end{theorem}
\begin{proof}
Since Banchoff's proof \cite{Ban70} of the critical point theorem for polyhedra is only
based on local identities that are trivially true also in $\Tt$, that proof holds with no
change for our case.
\par
The second part can be proved by considering that the genus $g$ of the surface $\oM$
contained inside a cube of side $R$ can be evaluated by reducing by homotopy the surface
to a graph and then evaluating the rank of the graph's first homology group.
The result follows from the consideration that every component contained in a inner
unitary cube contributes by $g$ to the total genus of the component contained in the cube
of radius $R$ and their number grows with $R^3$, while the cubes on the boundary provide a
smaller contribute but can be disregarded in the limit for $R\to\infty$ since their number
grows only as $R^2$.
\end{proof}
\subsection{Structure of foliations}
In the most general case, a constant 1-form $\omega$ induces on a triply periodic 
polyhedron $M$ both open and closed leaves. 
\par
Since being homotopic to zero is an open condition, leaves close enough to closed ones are
also closed; maximal components of closed leaves are always enclosed between a pair of
critical points of $\omega$ on $M$ and form either cylinders (when both critical points
are saddles) or discs (when one is a saddle and the other is a center) or spheres (when
they are both centers). The last two cases are topologically trivial: a disc covered by
closed leaves around a center is exactly a homotopy to a point of that component of $M$,
and if $M$ is a sphere then no open orbit can ever be induced on it by a closed
1-form. Hence, once the topologically trivial components are removed, what is left is a
collection of cylinders $\cC_i$ that separates a collection of subpolyhedra with boundary
$\cN_i$ filled by open leaves.
\begin{definition}
A {\sl genus-k component} of $\oM$ is a PL submanifold with boundary $\cN$ of $\oM$ such that:
\par
 \begin{enumerate}
   \begin{item}
   $\partial \cN$ is the finite disjoint union of plane parallel (topological) circles 
	homotopic to zero;
   \end{item}
   \begin{item}
   the closed polyhedron $\oN$, obtained by filling $\cN$'s holes with plane discs, 
	has genus k.
   \end{item}
 \end{enumerate}
\end{definition}
\begin{definition}
A {\sl Dynnikov decomposition} $\cZ$ of a polyhedron $\oM$ is a collection of subpolyhedra
with boundary $\{\cC_i,\cN_j\}$ of $\oM$ such that:
\par
 \begin{enumerate}
    \item every $\cC_i$ is homotopic either to a closed cylinder or to a closed disc;
    \item every $\cN_j$ is a genus-$k_j$ component of $\oM$;
    \item $i\neq j\implies\cN_i\cap\cN_j=\emptyset$, while every other pair of distinct subpolyhedra 
	of $\cZ$ that is not disjoint shares single a boundary component;
    \item $\oM=\bigcup_i\cC_i \bigcup_j\cN_j$.
 \end{enumerate}
The genus and rank of a Dynnikov decomposition $\cZ$ are, respectively, the highest genus
and rank of the $\oN_j$s contained in $\cZ$; a genus-1 rank-2 Dynnikov decomposition is
called ``Zorich decomposition''.
\end{definition}
%
From the considerations above it is clear that every constant 1-form $\omega$ in general
position with respect to $\oM$, i.e. such that there are no saddle connections, induces
naturally a Dynnikov decomposition $\cZ$ on $\oM$ where all $\cC_i$ are filled by the
closed leaves and all $\cN_j$ by the open ones.
%
\par
Under the same assumptions, to every focus it is associated a saddle that ``cancels'' it,
namely there's a homotopy of the surface that gets rid of the pair saddle-focus without
modifying the topology of the open orbits nearby; we call the saddle-type critical points
that are left ``topological saddles'', since it is them that contribute to the Euler
characteristic of the surface.
\par
In the smooth generic case, the number of $\Gamma$-inequivalent topological saddles,
i.e. the number of topological saddles in $\oM$, is of course exactly equal to $2g-2$; in
the PL case however, since we can have generic saddles with higher multiplicity, the
Euler characteristic is only an upper bound for the number of topological saddles, as
$\chi(\oM)=2-2g=\sum_{k\in\bZ}k\cdot\#\{[p]|i(\omega,p)=k\}$.
\subsection{The close-to-rational case}
\begin{definition}
By ``irrationality degree'' of a closed 1-form $\Omega$ on a piecewise smooth manifold $M$
we mean the number of rationally independent integrals of $\Omega$ over any base of the
homology integer 1-cycles of $M$: $\irr(\Omega)=\dim_\bQ\span_\bQ\{\int_\gamma\Omega\}$,
$\gamma\in H_1(M,\bZ)$; 1-irrational forms are also called ``rational''.
\end{definition}
The structure of foliations induced on a periodic polyhedron by rational 1-forms is much simpler 
than in the generic case since in this case all leaves are periodic; moreover, the appearance
of a topological invariant enforces all 1-forms close enough to rational to show the same
behaviour.
\par
Rational 1-forms are rather special because foliations induced on $\oM$ by them are also
induced by well-defined circle-valued functions on the polyhedron: indeed, if
$\{\gamma_i\}_{i=1,2,3}$ is a base for $H_1(\Tt,\bR)$ and $N$ an integer big enough
such that $N\int_{\gamma_i}\omega\in\bZ$, $i=1,2,3$, and $p_0\in\oM$, then the function
$$
\vbox{
  \halign{
	$#$\hfill&\hskip 2pt\hfill$#$\hfill&\hskip 2pt$#$\hfill&\hskip 2pt\hfill$#$\hfill\cr
	f:&M&\longrightarrow&\bS^1\cr
	&p&\mapsto&\exp(iN\displaystyle\int_{p_0}^p \omega)\cr
  }
}
$$ 
is well defined and differentiable, and its differential $df_p=iN\omega_p f(p)$ is
proportional to $\omega$; since $f$ is never zero, the set $df=0$ on $M$ coincides with
the restriction to $M$ of $\omega=0$.
\par
\begin{figure}
\begin{center}
\includegraphics[bb=0 300 612 692,width=6cm]{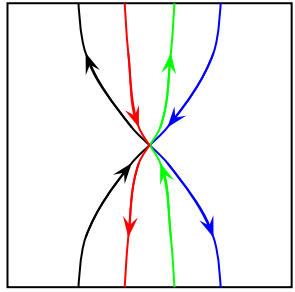}
\end{center}
\caption{%
A (hypothetical) saddle point where four open leaves meet simultaneously: the two branches
of each of the inner open leaves cannot be adjacent to each other, so this picture cannot
come from the section of a locally euclidean surface.
}
\label{fig:4ol}
\end{figure}

This shows that all leaves induced by rational 1-forms on a polyhedron $\oM\subset\Tt$ are
compact and therefore corresponding leaves on $M\subset\Rt$ will be open or closed
according to their homology class: leaves homotopic to zero in $\Tt$ will remain so in any
covering, all others will open up in $\Tt$'s universal covering.
\begin{definition}
From now on we will refer to closed leaves non homotopic to zero as ``periodic'' or, more
generically, ``open'' leaves, so that by ``closed leaf'' we will implicitly mean a closed
leaf homotopic to zero.
\end{definition}
\begin{lemma}
\label{lemma:4ol}
Be $\omega$ a rational 1-form in general position with respect to a polyhedron $\oM$: then
no more than two open leaves can collide at any saddle point.
\end{lemma}
\begin{proof}
The 1-form $\omega$ foliates $\Tt$ in a 1-parameter family of embedded 2-tori, so that all
open leaves at the same level of $\omega$ are parallel (i.e. they represent the same
1-cycle modulo sign); moreover the number of open leaves on every level is even, since
they are all indivisible and their sum must be zero.
\par
It is easy to check, just by drawing pictures, that it is possible to have saddles of any
index with closed leaves, and adding a single pair of open leaves does not change this
situation; their presence though does not allow the presence of any other pair, since 
in any saddle point the two extremes of the same (critical) leaf must appear next to
each other and this is of course impossible for all open leaves, apart for the two most 
external ones (e.g. see Figure~\ref{fig:4ol}).
\end{proof}
\begin{theorem}
\label{thm:Zor}
The Dynnikov decomposition induced on $\oM$ by any constant rational 1-form in general
position is Zorich.
\end{theorem}
\begin{proof}
Since all (non critical) leaves induced by $\omega$ on $\oM$ are circles, the critical
points of $\omega$ on $\oM$ determine a subdivision of the polyhedron in the connected sum
of a finite number of cobordisms. In the smooth case, the only non trivial cobordisms are
elementary (i.e. pants-like) and Zorich \cite{Zor84} proved that in the boundary of each
pant there is at least one of the three boundary loops that is homotopic to zero (in
$\Tt$), from which the theorem follows easily.
\par
In the polyhedra case the cobordisms are not necessarily elementary, i.e. more than two
leaves may collide at a saddle point, but lemma~\ref{lemma:4ol} grants us that even in
this case no more than two loops in the cobordism may be open.
This shows that the components of open leaves have genus 1 and therefore the decomposition
induced on $\oM$ by $\omega$ is Zorich.
\end{proof}
The presence of a genus-1 rank-2 component of open leaves is a very strong condition and leads
to the appearance of a topological quantity associated to the foliation~\cite{Dyn97}:
\begin{theorem}
A Dynnikov decomposition $\cZ$ of $\oM$ with rank 2 has genus 1 if and only if at least one of
its $\cN_j$ embedded with rank 2 has genus 1. In this case, all $\cN_j$ are genus-1
components of $\oM$ and those embedded with rank 2 represent, modulo sign, the same
indivisible non-zero 2-homology class in $\Tt$.
\end{theorem}
\begin{definition}
We call ``soul'' of a Zorich decomposition $\cZ$ the ``unsigned'' non-trivial
indivisible homology class $l\in PH_2(\Tt,\bZ)$ common to all rank-2 components of
$\cZ$.
\end{definition}
\begin{corollary}
Be $\omega$ a rational 1-form inducing on $\oM$ a rank-2 Zorich decomposition $\cZ$:
then any 1-form $\omega^\prime$ close enough to $\omega$ induces on $\oM$ a rank-2 Zorich
decomposition $\cZ^\prime$ and this decomposition is homotopic to $\cZ$.
In particular, all such decompositions share the same soul, i.e. the soul is a locally
constant function of the pair $(\oM,\omega)$. 
\end{corollary}
\begin{proof}
The leaves at the boundary between $\oM$'s genus-1 components of periodic leaves $N_i$ and
the cilinders of closed leaves are themselves closed and therefore stable under small
perturbations of $\omega$'s direction; consequentely, no cilinder of closed orbits will
disappear in a whole neighborhood of $\omega\in\RPt$. Since no two leaves of a foliation
can intersect, this means that open leaves are bounded to genus-1 components of $\oM$,
homotopic to the $\cZ$ ones, even for all 1-forms close enough to $\omega$.
\par
Finally, since the homology class $l$ of these rank-2 genus-1 components is integer and
they change continously, this $l$ must be the same (modulo sign) for all them.
\end{proof}
\par
The soul $l$ of a Zorich decomposition is a fundamental invariant since its knowledge is
enough to describe the asymptotic behaviour of the open leaves. Indeed, the fact that
$[\cN_i]=\pm l$ means that, in $\Rt$, the lift $\hat\cN_i$ of $\cN_i$ lies between a pair
of parallel planes perpendicular to $l$ (seen as a direction in $\Rt$), so that the lift
to $\Rt$ of an open leaf, namely an open intersection of $\hat\cN_i$ with a plane
perpendicular to $\omega$ (seen as a vector in $\Rt$), is a curve contained in a plane
strip of finite width; Dynnikov~\cite{Dyn92} showed how these conditions are enough to
conclude that the leaf is actually a finite deformation of a straight line whose direction
is the axis of the strip, namely ``$\omega\times l$''.
\par
Since the foliation induced by $\omega$ is determined just by its direction $\bH$ and the
soul itself can be interpreted as a direction in $\Rt$, once we fix a surface $M$ we can
think the soul application as a locally constant application $\soul_{\oM}:\RPt\to\RPt$.
We will show in next section that, for a generic polyhedron $\oM$, this map is well
defined on the whole projective plane except for a set of meaure zero (``ergodic
directions'') and its image amounts to a finite number of points.
\begin{definition}
The non-empty level sets $\cD_l(\oM)=\soul_{\oM}^{-1}(l)$ corresponding to non-zero values of $l$
are called ``islands'' or ``stability zones'' of $\oM$. The union of all stability zones
$\cS(\oM)=\cup_l\cD_l(\oM)$ is called the ``Stereographic Map'' (SM) of $\cM$.
\end{definition}
Strange as it seems, this whole construction is the natural model for a very concrete
physical phenomenon, namely the magnetoresistance in normal metals, with the periodic
surface being the ``Fermi Surface'' of a metal, the 1-form a strong constant homogenous
external magnetic field and the Fermi Surface leaves the orbits of the momenta of the metal's
quasi-electrons.  The presence of open leaves is detectable experimentally, and so is each
stability zone (provided at least two points of it are measured \cite{NM98}). 
Experimental plots of the SM have been produced in Sixties and Seventies for about thirty
metals but only recently the first two SM, relative to the Fermi Surfaces of Au and Ag,
have been reproduced theoretically from first first principles~\cite{DL04,DL05a}.
\subsection{The generic (3-irrational) case}
In the previous section we showed that the structure of the foliation in a neighborhood of
rational directions is rather simple: either all leaves are closed or their lift is
strongly asymptotic to a straight line~\cite{Dyn97}, i.e. it is contained in a
finite-width plane strip and crosses it from one end to the other.
Still, this is far from being enough, since the sole density of an open set does not
preclude its measure from being small.
\par
We will prove in this section that the one above is nevertheless the generic situation,
namely:
\begin{theorem}
\label{thm:Dyn}
The set of directions in $\RPt$ inducing Dynnikov decompositions of genus bigger than one
on a generic polyhedron $\oM$ has measure zero.
\end{theorem}
To prove this theorem, Dynnikov studied the structure of the foliation induced by a
3-irrational 1-form $\omega$ on 1-parameter family of surfaces $M_e=f^{-1}(e)$, where $f$ is
a triply periodic Morse function such that for almost all values of $f$ no more than one
critical point of $\omega$ lies on the same leaf. Here, we will repeat Dynnikov's steps
assuming $f$ to be a generic triply periodic PL function, so that almost all of its level
surfaces are embedded polyhedra, satisfying the same genericity condition with respect to
$\omega$.
\par
The following two lemmas by Dynnikov extend with no change to the polyhedra case:
\begin{lemma}
\label{thm:interval}
Given a triply periodic PL function, the set of values for which a constant 1-form $\omega$
induces open leaves filling rank-2 components is a closed connected interval
$[e_1(\omega),e_2(\omega)]$. The functions $e$ and $E$ are continous on the whole $\RPd$.
\end{lemma}
\begin{definition}
The Dynnikov index $w$ of a critical point $c$ of $\omega=h^\alpha dp_\alpha$ with respect
to an oriented polyhedron $\oM$ is the product of the ``Hamiltonian'' index of the critical
point ($+1$ for centers and $-1$ for saddles) times the sign of the scalar product between
$\bH=(h^\alpha)$ and any of the normals to the faces adjacent to the critical vertex.
\end{definition}
\begin{lemma}
The curve $\gamma_{\omega,f}(e)=\sum_i w_i c_i(e)$ is a well defined loop in $\Tt$ and the
quantity 
$$\overline{\chi}_{\omega,f}(e)=\int_{-\infty}^e i^*_{\gamma_{\omega,f}}\omega$$ 
is equal to the density of closed leaves on any leaf induced by $\omega$ on $\Tt$.
\par
Moreover, if $\{h^+_j\}$ (resp. $\{h^-_k\}$) are the heigths of the ``positive''
(resp. ''negative'') cylinders of closed leaves on $\oM$, namely those that contain points
with smaller (resp. bigger) values of $f$, it turns out that it is possible to choose
$\Rt$ representative $\hat c_i$ of the critical points so that
$$\overline{\chi}_{\omega,f}(e)=\sum_{i} <\bH, w_i \hat c_i(e)>=\sum h^+_j - \sum h^-_k$$
\end{lemma}
In the weighted sum of all critical points are also contained the pairs center-saddle,
that can be removed by homotopy and have nothing to do with the topology of the foliation.
If we do not include them in the summation, we are left with a new quantity that tells us
the density of non-trivial (in $\oM$) closed leaves and therefore is able to spot whether
there are cylinders of non-zero heights. In particular, to have a full ergodic situation,
i.e. a leaf dense on the whole $\oM$, this function must necessarily be zero.
\begin{lemma}
Be $\{c^{top}_i\}$ the subset of the topological saddles of the pair $(\omega,f)$. Then
the reduced curve $\gamma^{top}_{\omega,f}(e)=\sum_i w_i c^{top}_i(e)$ is also
well-defined in $\Tt$ and the reduced Euler density
$\overline{\chi}^{top}_{\omega,f}(e)=\int_{-\infty}^e i^*_{\gamma^{top}_{\omega,f}}\omega$
is equal to the density of closed leaves non homotopic to zero in $\oM$.
\par
The function $\overline{\chi}^{top}_{\omega,f}$ is a strictly increasing non continous PL
function with respect to $e$; its discontinuity points are exactly the non-proper values
of $f$.
\end{lemma}
\begin{proof}
The original curve $\gamma_{\omega,f}(e)$ is well-defined because $\sum
w_i=0$~\cite{Dyn97}; since the Dynnikov indices of a pair saddle-center are opposite, the
restricted sum $\sum_{top} w_i$ is still zero and therefore also
$\gamma^{top}_{\omega,f}(e)$ is well defined.
\par
Considering only the sum of the topological saddles is equivalent to cancel from a leaf
all closed leaves that are homotopic to zero in $\oM$, so that the averaged Euler
characteristic now is only relative to the non-trivial closed leaves only.
\par
Now consider a positive cylinder: since the values inside are smaller with respect to the
values on the cylinder, the height of the cylinder increases with $e$; for the same
reason, negative cylinders decrease their height.
Since $\overline{\chi}^{top}_{\omega,f}(e)=\sum_{top} h^+_j - \sum_{top} h^-_k$, then the
function $\overline{\chi}^{top}_{\omega,f}(e)$ is strictly increasing in its continuity
points.  The function fails to be continous when new pairs of topological saddles are
created or destroyed, namely in the non-proper values of $f$, and it jumps exactly by the
height of the cylinder(s) created or destroyed.
\end{proof}
\begin{corollary}
If $\omega$ induces on $\oM_e$ a Dynnikov decomposition of genus bigger than 1, then
$e_1(\omega)=e_2(\omega)=e_0$, i.e. $\omega$ induces only closed leaves at any other
level.
\end{corollary}
\begin{corollary}
At almost all levels of $f$ the measure of the ``ergodic'' directions is zero.
\end{corollary}
We have now all ingredients to prove theorem~\ref{thm:Dyn}:
\begin{proof}
We will just discuss the case of ``full ergodicity'', namely the case of directions giving
rise to leaves dense on the whole polyhedron; the same line of arguments extends to the
lesser ergodic cases.
\par
If $\omega$ induces on $\oM$ fully ergodic leaves, then of course no cylinder can appear
and therefore, for any Morse PL function $f$ such that $\oM=f^{-1}(0)$, it must happen that
$\overline{\chi}^{top}_{\omega,f}(0)=0$. Let us consider now $\overline{\chi}^{top}$ as a
function of $\omega$ and $e$: then the surface
$X=(\overline{\chi}^{top}_f)^{-1}(0)\subset\RPt\times\bR$, for a generic function $f$, is
transversal to the sections $\RPt\times\{e\}$: indeed in a projective chart, say $h^z=1$,
we have that $\partial_{h^x}\overline{\chi}^{top}_{\omega,f}=\sum_i w_i (c^{top}_i)^x$ and
similarly for $h^y$, so that the points on $X$ where the gradient is zero are exactly the
points where $\gamma^{top}_{\omega,f}(e)=\sum_i w_i c^{top}_i=0$.
This condition is non-generic, that finally proves the claim of the theorem.
\end{proof}
\subsection{Structure of the Stereographic Map of surfaces and triply periodic functions}
The results of the previous two sections show that the SM $\cS(\oM)$ of a generic surface
is the disjoint union of a countable set of open sets $\cD_l(\oM)$ (``islands''), each
labeled by a $l\in PH_2(\Tt,\bZ)$, immersed in a sea of directions that give rise only to
closed leaves. According to our intuition of the system and the numerical experiments made
to date, we conjecture that generically the number of islands is finite, but no rigorous
proof of this fact exists. As matter of fact, this structure is exactly the one guessed,
from symmetry consideration, by the physicist I.M. Lifschitz and his Karkov school about
fifty years ago~\cite{LP59,LP60}.
\par
The boundaries of the islands are reached when the last pair of genus-1 components collide
because of a cylinder collapse and therefore are charaterized by the presence of (at
least) a pair of inequivalent critical points on the same leaf. The set of these
directions is the countable union of the curves $<\bH,\hat c_i-\hat c_j>=0$, $i\neq j$; such
curves in the polyhedra case are all straight lines, so that every island is actually
a (not necessarily convex) polygon. There are reasons to believe that these polygons
are convex for low genus, i.e. at least for genus 3 and 4, but it is easy to build examples
of high-genus polyhedra with islands that are either non connected or connected but non convex
or even connected but non simply connected.
\par
A crucial observation by Dynnikov~\cite{Dyn97} allows to associate a SM also to triply
periodic functions:
\begin{theorem}
Be $f$ a Morse triply periodic function: then, if $\omega$ induces on $\oM_{e0}$ a Zorich
decomposition $\cZ_{e0}$, it induces a Zorich decomposition $\cZ_e$ for all
$e\in(e_1(\omega),e_2(\omega))$ and all these decomposition share the same soul.
\end{theorem}
\begin{definition}
The island $\cD_l(f)$ corresponding to the label $l$ is the union of the corresponding
islands of $f$'s level sets: $\cD_l(f)=\cup_{e\in\bR}\cD_l(\oM_e)$.
The SM of $f$ is the union of all its islands: $\cS(f)=\cup_l\cD_l(f)$
\end{definition}
The SM corresponding to functions are generically drammatically different from those
corresponding to surfaces. 
\par
First of all, since rational directions induce necessarily Zorich decompositions, the set
of islands $\cS(f)$ is now always dense in $\RPt$, and of course
lemma~\ref{thm:interval} also tells us that in this case the sea of directions giving rise
to closed leaves only dried up, since every direction either belongs to an islands (or its
boundary) or is ergodic.
\par
Moreover, the following property shows that the islands can be sorted in a rather complex
way:
\begin{theorem}
Generically every two zones meet transversally and in a countable number of points.
\end{theorem}
\begin{proof}
No point belonging at the same time to two different zones $\cD_{l_1}$ and $\cD_{l_2}$ can
have irrationality degree bigger than 2, since it must contain the integer direction
$l_1\times l_2$. In the smooth case this would be enough, since the boundaries are smooth
curves and they generically contain only a countable number of 2-irrational points and no
rational point. In the PL case though the boundaries are actually segments of straight
lines and threfore they are actually contained in the set of the directions with
irrationality degree smaller than one.
\par
Nevertheless, the theorem holds also in this case for the following reason: since every
direction at the boundary between two zones is perpendicular to the direction $l_1\times
l_2$, then their set is the straight line (in $\RPt$) passing through $l_1$ and
$l_2$. Generically none of the two labels falls on the boundary and therefore two zones
can meet in a number of points not bigger than the number of sides of the island with the
smaller number of sides.
\end{proof}
\begin{corollary}
Either there is a single zone, i.e. it exists a label $l$ such that $\cD_l(f)=\RPt$, or there
are countably many zones and they are dense in the whole projective plane.
\end{corollary}
Since the islands meet trasversally, the non trivial SM will look like 2-dimensional
Cantor sets.  Retrieving numerically such fractals is not trivial since it involves, in
general, to analyze the system at several values of $f$, but there is a class of
interesting (non generic) cases in which it is actually enough to analyze a single level
surface of $f$ to get the entire fractal picture:
\begin{theorem}
\label{thm:intext}
Be $\oM$ a polyedron whose interior is equal to its exterior, modulo the group
$G\simeq\Rt\times\bZ_2$ of translations and inversion of the three axes.  Then
$\cS(\oM)=\cS(f)$ for any function having $\oM$ as (connected component of a) level set.
\end{theorem}
\begin{proof}
By symmetry, we can build a function $f$ such that $\oM=f^{-1}(0)$ and that $f^{-1}(e)$ is
equal to $f^{-1}(-e)$, modulo $G$. Since also the bundles of parallel planes are invariant
by $G$, it turns out that for such an $f$ the interval $I(\omega)$ of existence of open
orbits relative to any 1-form $\omega$ is of the form $I(\omega)=[-a(\omega),a(\omega)]$
and therefore $0\in I(\omega)$, $\forall\omega\in\RPt$, i.e. at the zero level every
$\omega$ induces open leaves.
\end{proof}
In particular, all triply periodic functions $f$ such that $f(c-x,c-y,c-z)=-f(x,y,z)$
belong to this class, and indeed the only fractals analyzed numerically to date are
relative to this kind of functions.
\par
Finally, we cite an important property that ties the set of all labels relative to the islands
of the SM of a function with the set of ergodic direction~\cite{DL03b,DL05b}:
\begin{theorem}
\label{thm:labels}
The closure of the set of all labels is the disjoint union of the
set of all zones boundaries and the set of ergodic directions.
\end{theorem}
\section{A concrete case study}
\label{sec:sp}
As pointed out in the previous section, to date a picture of the fractal has been
numerically produced for only two functions, an analytical one and a piece-wise quadratic
one~\cite{DL03b}.
\par
The analytical one is $\mf(x,y,z)=\cos(2\pi x)+\cos(2\pi y)+\cos(2\pi z)$, invariant with
respect to translations by integers $\Gamma=\bZ^3\subset\Rt$, that gives rise
to genus-3 level surfaces in the range $(-1,1)$ and spheres at every other non-critical
level.  This function represents the simplest non-trivial case possible from the
topological point of view, since any triply periodic connected surface of genus smaller
than 3 lies between two parallel planes and therefore the aymptotics of plane sections is
easily found.  Its zero level is rather special: it is known as the Schwarz primitive
function (or plumber's nightmare) and was studied by Schwarz in 1890 as one of the first
examples of triply periodic minimal surface.
\par
From the computational point of view, $\cS\cP=\mf^{-1}(0)$ has three important properties:
\begin{enumerate}
\item[{\textbf SP1}] its interior is a translate of its exterior, so that $\cS(\cS\cP)=\cS(\mf)$;
\item[{\textbf SP2}] it is invariant with respect to the natural action of the tetrahedral group $T_d$ 
      on the unitary cube, so that the whole SM can be obtained, for example, extending 
      by symmetry to the whole $\RPt$ the data obtained for 
      the triangle with vertices $[(0,0,1)]$, $[(1,0,1)]$ and $[(1,1,1)]$;
\item[{\textbf SP3}] the two cylinders, one negative and one positive, have the same height, so that it
      is enough to examine just one of the four topological critical points at the base
      of the two cylinders to retrieve all information about the structure of the foliation.
\end{enumerate}
\par
The piecewise quadratic function is $\mg(x,y,z)=\sum\bar \mg(x^i)$, where $\bar\mg$ is the
simplest piecewise quadratic function having the same symmetries of the cosine
function. Its level sets have the same behaviour as the function above but the expression
of the critical points as function of the direction of the 1-form and the level of the
function is so simple that allows a comparison between analytical and numerical data also
at levels different from zero.
\par
Nevertheless, in both cases the number of triangles needed to describe in sufficient
detail the surface is so big (between $10^5$ and $10^6$) to make impossible to improve the
resolution of the results obtained in~\cite{DL03b}, at least until a new algorithm is found
or some ``ad hoc'' trick used.
\subsection{The polyhedron}
\begin{figure}
\begin{center}
\includegraphics[bb=200 512 400 712,width=3.8cm]{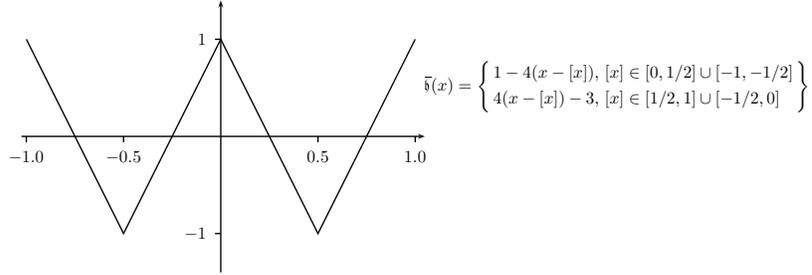}
\end{center}
\caption{%
The simplest PL approximation of the cosine function.
}
\label{fig:lincos}
\end{figure}
A natural way to improve the resolution of the numerical analysis of the problem is to
consider PL functions, since in this case the number of triangles needed to describe their
level surfaces can be as low as of the order of $10^1$. Moreover, the description of the
surfaces is in this case exact rather than an approximation.
\par
\begin{figure}
\begin{center}
\includegraphics[width=5cm]{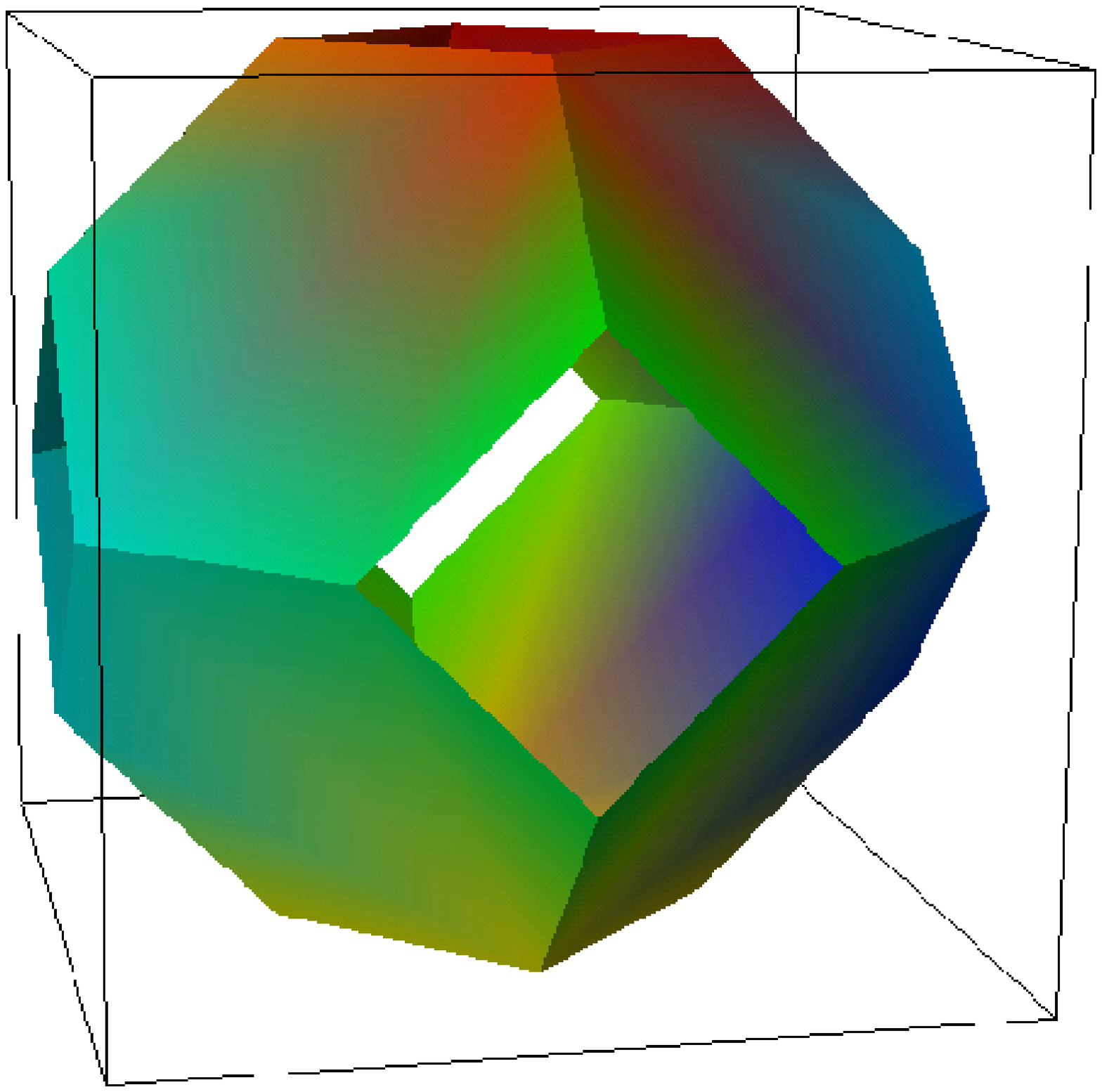}
\hskip1.cm
\includegraphics[width=5cm]{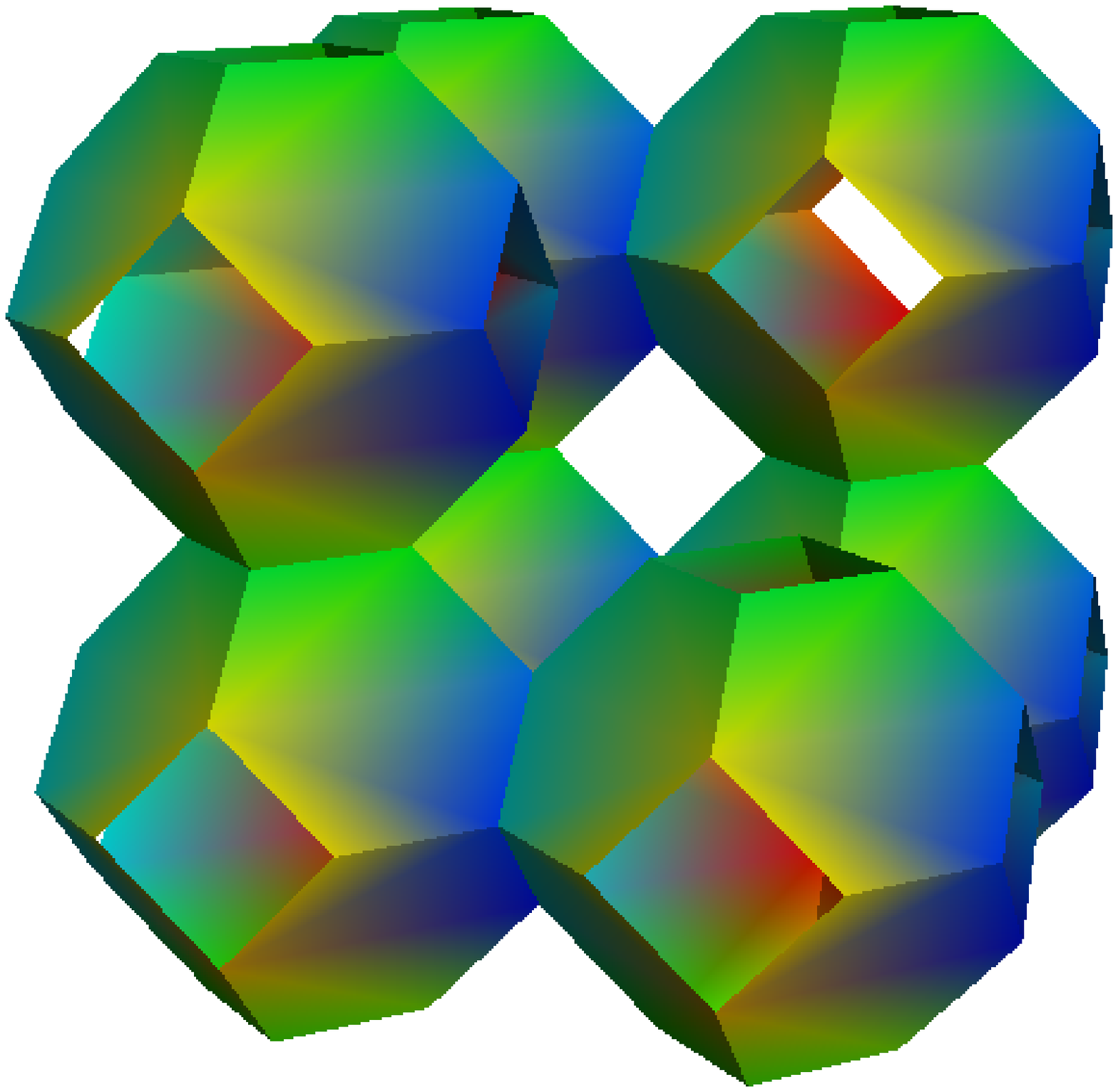}
\end{center}
\caption{%
Plot of the Truncated Octahedron inside $\Tt$ (left) and of part of its image in the
universal covering.
}
\label{fig:pol}
\end{figure}

The simplest case to study, in order to take advantage of the extremely convenient
properties evidenced in the $\cS\cP$ case, is the PL function
$\mh(x,y,z)=\sum_i\bar\mh(x^i)$, where $\bar\mh$ is the function shown in
Fig.~\ref{fig:lincos}.
The polyhedron $\cP_0=\mh^{-1}(0)$ is a PL embedding of a genus-3 surface in $\Tt$ having
8 hexagonal faces, 20 edges and 12 vertices; the smallest triangulation for $\cP_0$ takes
$4\cdot8=32$ triangles, $32\cdot3/2=48$ edges and $12$ vertices, that gives the expected
Euler characteristic $\chi(\cP_0)=F-E+V=-4=2-2\cdot3$.
The basic cell of the lift $\overline\cP_0\subset\Rt$ in the unit cube is a truncated octahedron
(Fig.~\ref{fig:pol}); it is noteworthy to notice that, like its smooth analog, also this surface
is, in the discrete sense, a minimal surface~\cite{Way04}.
\par
Since exactly four edges meet at every vertex, given any 1-form $\omega$ in general
position with respect to $\cP_0$ only saddles with index $-1$ may arise and therefore
there are always four vertices that are critical for such $\omega$. If the point $p_1$ is
such a vertex, then the remaining three critical points are $p_{2,3}=(1/2,1/2,1/2)\pm p_1$
and $p_4=(1,1,1)-p_1$. In particular in our numerical study we sample the set of 1-forms
$\omega=(h^x,h^y,h^z)$ such that $h^x/h^z\in[0,1]$ and $h^y/h^z\in[0,1]$, for which the four
critical points are $(0,.5,.75)$, $(.5,1,1.25)$, $(.5,0,.25)$ and $(1,.5,.25)$.
\subsection{The algorithm}
In order to generate an approximate picture of the fractal, it is enough to produce an
algorithm able to evaluate the label, if any, associated to a given 1-form. Since
obviously no calculator can deal with irrational numbers, the numerical study will be
limited to rational 1-forms; luckily this is not a big restriction, since anyway rational
directions are dense in every stability zone.
\par
Note that the algorithm used for the numerical study of the PL case is a simplified
version of the more general algorithm we developed to study smooth surfaces of genus
three~\cite{DL03b}, since in this case we know a priori the position of all critical points
and moreover we know their position exactly, so that we do not have to correct ``by hand''
the topology of the critical section.
\par
The basic idea to retrieve the label, as suggested to me by I. Dynnikov, is that the soul
$l\in PH_2(\Tt,\bZ)$ associated to the Zorich decomposition $\cZ$ induced by $\omega\in\cD_l$ is
in 1-1 correspondance with the rank-2 sublattice of $H_1(\Tt,\bZ)$ obtained as the image,
through the map $i_*:H_1(\oM,\bZ)\to H_1(\Tt,\bZ)$, of the open leaves in $\oM$ that have
zero intersection number with the closed leaves populating the cylinders of $\cZ$.
Indeed, every cycle lying on the interior of a $N_j$ component of $\cZ$ has no
intersection with the closed leaves that form the cylinders $\cC_i$, and since all $N_j$
are homologous to each other (modulo sign) the image of all these cycles in $\Tt$ must
have rank-2; on the other side, there is an obvious 1-1 correspondance between rank-2
sublattices of $H_1(\Tt,\bZ)=\Gamma\simeq\Zt$ and the 2-tori embedded in $\Tt=\Rt/\Gamma$,
since every such 2-torus can be spanned by a pair of independent rational directions and
viceversa.
\par
The following algorithm {\textbf N} works for genus-3 polyhedra satisfying properties
{\textbf SP1-3}, in particular for $\cP_0$:
\par
{\textbf Input}: $\oM$ - the polyhedron; $\omega=(l,m,n)\in\Zt$ - the 1-form; $x$ - a critical
point of $\omega$ with respect to $\oM$; $\pi_{\omega,x}$ - the plane perpendicular to
$(l,m,n)$ and passing through $x$.
\par
{\textbf Output}: $c_{1,2}$ - the two critical loops%
\footnote{since numerically we can study only the 1-rational case, in $\Tt$ the saddles are
always wedges of circles, i.e. all critical branches close back to the critical point; according to
whether these loops are or not homotopic to zero in $\Tt$, their $\Rt$ lift will be open or close.}%
; $h_{1,2}$ - the homology classes of $c_{1,2}$ in $\oM$; $H_{1,2}$ - the homology classes
of $c_{1,2}$ in $\Tt$.
\par
{\textbf Algorithm}:
\begin{itemize}
\item[{\textbf N1}] retrieve the intersection between $M$ and $\pi_{\omega,x}$;
\item[{\textbf N2}] check that there are exactly four critical branches and follow them
	by periodicity, otherwise exit;
\item[{\textbf N3}] if no other critical point is met along the path, so that the four
	branches are arranged in a pair of critical loops, store the two loops in the 
	variables $c_{1,2}$, otherwise exit;
\item[{\textbf N4}] evaluate the homology class of $c_{1,2}$ in $\Tt$ and in $M$
	(this is actually done while executing N2 to speed up the computations time);
\item[{\textbf N5}] if the saddle is half-open, i.e. if exactly one among $H_{1,2}$ is zero,
	then associate to $\omega$ the complementary $h$ triple, otherwise exit.
\end{itemize}
The main outcome of the algorithm is of course the label associated to $\omega$. The fact that
this label is a triple of integers is very important, since an integer evaluated numerically
with an error smaller than .5 becomes actually an exact measure.
\par
We implemented this algorithm in a C++ library named NTC \footnote{http://ntc.sf.net/} built
over an Open Source C++ library named VTK \footnote{http://www.vtk.org/}. The choice of the 
language comes from the fact that VTK provides the basic geometric environment and algorithms 
needed by the problem, mainly the capability of generating meshes for isosurfaces and
evaluating intersections between geometric objects. The inheritance mechanism of the 
C++ language allows to use transparently all functions of a library, hence we used VTK as
a starting point and implemented in NTC the routines to deal with periodicity and
evaluate the homology classes.
\par
No serious attempt to evaluate the error on such calculations has been made to date since
no need for it manifested. As check for the reliablity of the result were rather used indirect
evidences:
\begin{itemize}
\item the agreement of the biggest zones with their analytical boundary
(Fig.~\ref{fig:sm100}b), obtained through the independent algorithm {\textbf A};
\item the symmetry of the final picture with respect to the diagonal
(Fig.~\ref{fig:sm1000}), symmetry that was in no way used in the numerical calculations;
\item the agreement of the fractal picture with the labels plot (Fig.~\ref{fig:labels}).
\end{itemize}
\par	
The exploration of the SM was performed in the square $[0,1]^2$ of the projective chart
$(\omega_x/\omega_z,\omega_y/\omega_z)$ by evaluating the label associated to every
direction at the vertex of a uniform grid of step $r$ and it was repeated for the values
$r=10^2,10^3,10^4$.  Samplings with $r=10^2,10^3$ have been succesfully performed also for
the previous two functions \cite{DL03b} but the CPU time needed to reach $r=10^4$ in that
case was way too big.  It is because of the rather small number of triangles needed to
describe $\cP_0$ that the computation became doable.
%
\subsection{Numerical results for $r=10^2$}
\begin{figure}
\begin{center}
\includegraphics[width=12cm]{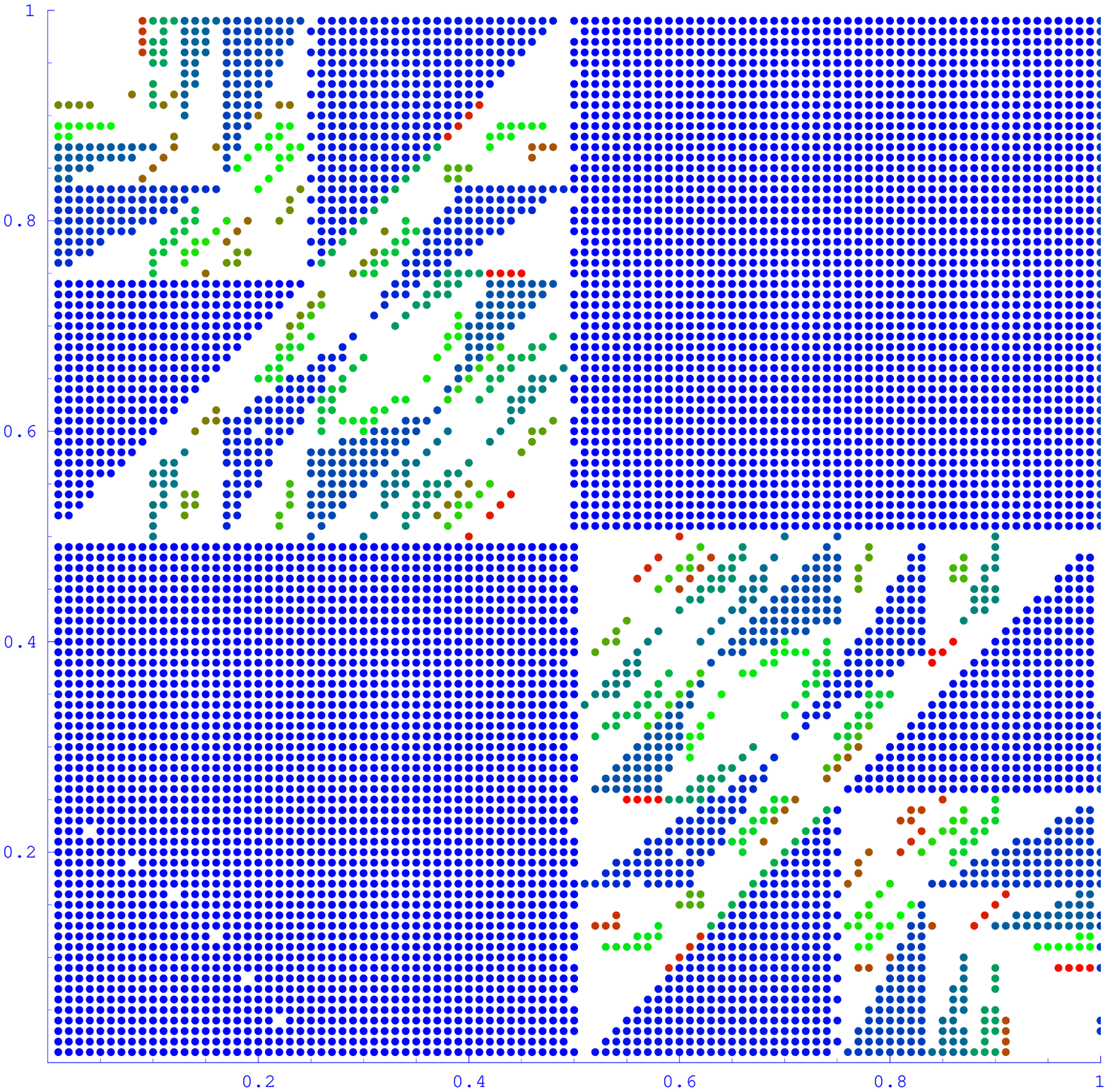}
\end{center}
\caption{%
Numerical plot of the square $[0,1]^2$ of the SM $\cS(\cP_0)=\cS(\mh)$ in the projective
chart $h^z=1$ at a resolution $r=10^2$.  The color of the islands goes from blue to red as
the norm of the label grows. In the picture are displayed the 106 islands with at least
four points out of the total 1741 islands found. The missing points that is possible to
see in the interior of some of the islands are due to failure of the numerical algorithm
{\textbf N}, e.g. because of the presence of saddle connections.  The running time for the
$10^4$ steps cycle needed to retrieve these data takes about $1h$ on a Pentium \~1GHz CPU.
}
\label{fig:sm100}
\end{figure}

This resolution is the lowest one that allows to have a hint of the structure of the
fractal.  About $r$ sections are needed to follow the critical branches for a generic
direction $(m,n,r)$, $m,n\in\{1,\dots,r\}$, that takes a time of $~.5s$ on a $~1GHz$ CPU
for the evaluaton of a single label and $10^4\times.5s\simeq1h$ for sampling the $10^4$
directions of the grid (Fig.~\ref{fig:sm100}).
\par
Even from this rough picture it is rather evident a further symmetry of the picture,
namely the one with respect to the antidiagonal of the square. This symmetry does not
come, like the others, from the tetrahedral group $T_d$ but it is rather of topological
nature. The numerical evidence is that, if a 1-form $(m,n,r)$ is labeled by $L$, then its
symmetric $(1-n,1-m,r)$ is labeled by $L+(1,1,0)$ but no proof of this fact is known.
\par
In order to verify the correctness of the algorithm, we found ``by hand'' the analytical
boundaries for the biggest zones and compared them with the numerical results
(Fig.~\ref{fig:sm100}b). 
The following algorithm, aimed at the cases similar to the cosine one, is a slight modification
of the original algorithm introducted and used by Dynnikov in~\cite{Dyn96}:
\par
{\textbf Algorithm A}
\begin{itemize}
\item[A1] fix a 1-form $\omega=(m,n,r)\in\bZ^3$ inside some zone 
	(e.g. extracting it from the 
	experimental, guessing it from symmetry arguments or simply by trial and error);
\item[A2] retrieve the critical section of $\omega$ passing through one of its
	critical points $c$ and make sure it is half open, otherwise exit;
\item[A3] evaluate the homology class $l$ of the closed critical leaf $C$;
\item[A4] rotate $\omega$ around some direction till the cylinder of which $C$
	is a base collapses and identify the critical point $c^\prime$ that is now
	connected to $c$ through a saddle connection;
\item[A5] the equation $<\bH,c-c^\prime>=0$ contains one of the sides of the island;
	follow it in one direction till four critical points fall over the critical closed
	leaf: this is the point when a sides and a new one start; repeat this
	step till the island boundary close up on themselves. 
\end{itemize}
\begin{figure}
\begin{center}
\includegraphics[bb=200 462 400 722,width=4.5cm]{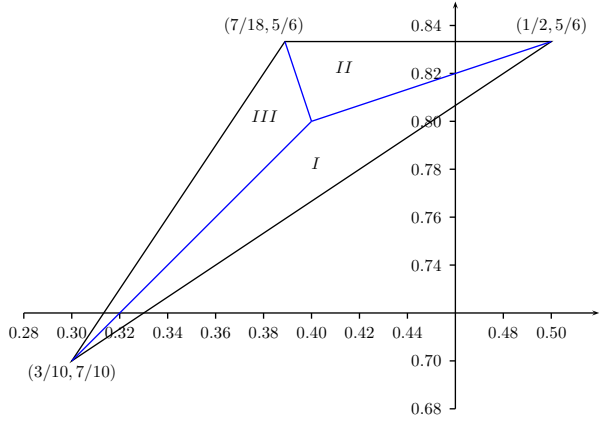}
\end{center}
\caption{%
A close-up of the island $\cD_{(2,4,5)}(\cP_0)=\cD_{(2,4,5)}(\mh)$.
Inside the island the pairs of critical points at the base of each cylinder are locally
constant. In correspondance to each sides, there are three different pairings sorted in
open subsets, labeled in the picture by roman numerals, separated by straight lines
segments correspoding to directions $\omega$ for which the bases of the positive and
negative cylinders collide, resulting in a saddle connection between two critical
points. These three segments meet in the single point $(.4,.8)$, that in this case happens
to be exactly the direction of the label (this property though is not generic).
In Figure~\ref{fig:sections} we show in detail the transition between two different pairs
within this island.
}
\label{fig:2.4.5}
\end{figure}
\par
\begin{figure}
\begin{center}
\vbox{\halign{\hfill#\hfill&\hfill#\hfill&\hfill#\hfill\cr
\includegraphics[width=4cm]{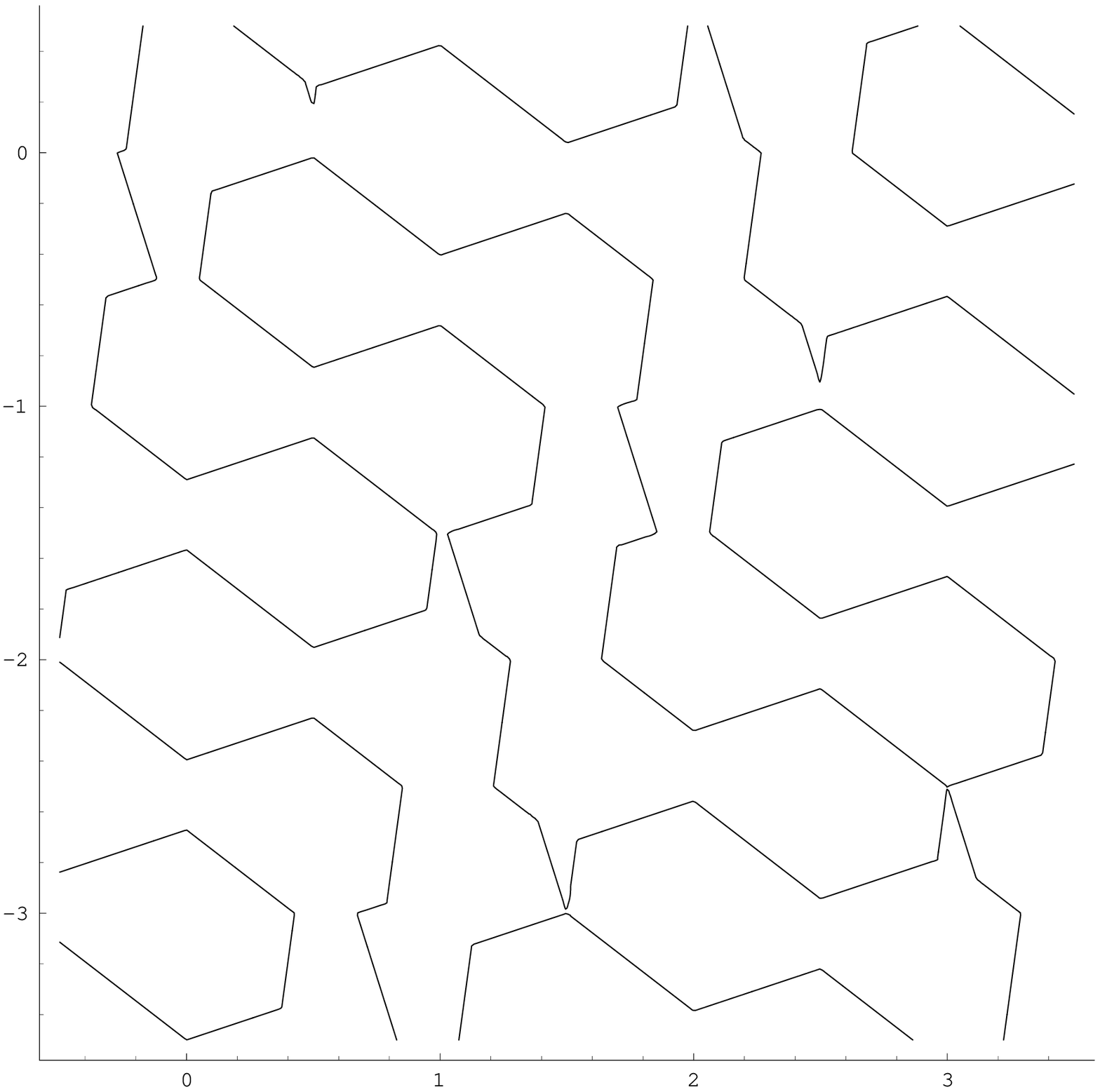}
&
\includegraphics[width=4cm]{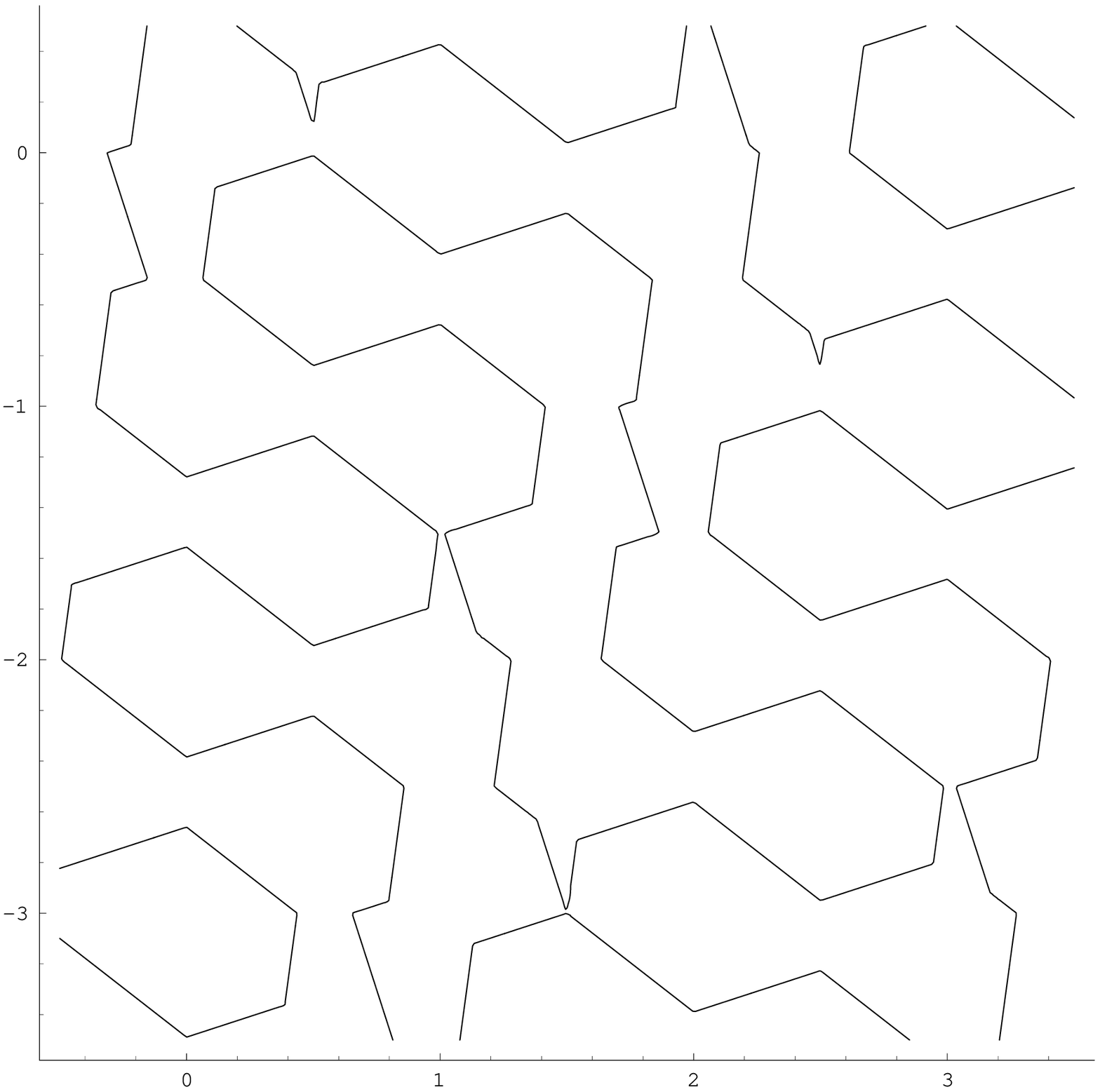}
&
\includegraphics[width=4cm]{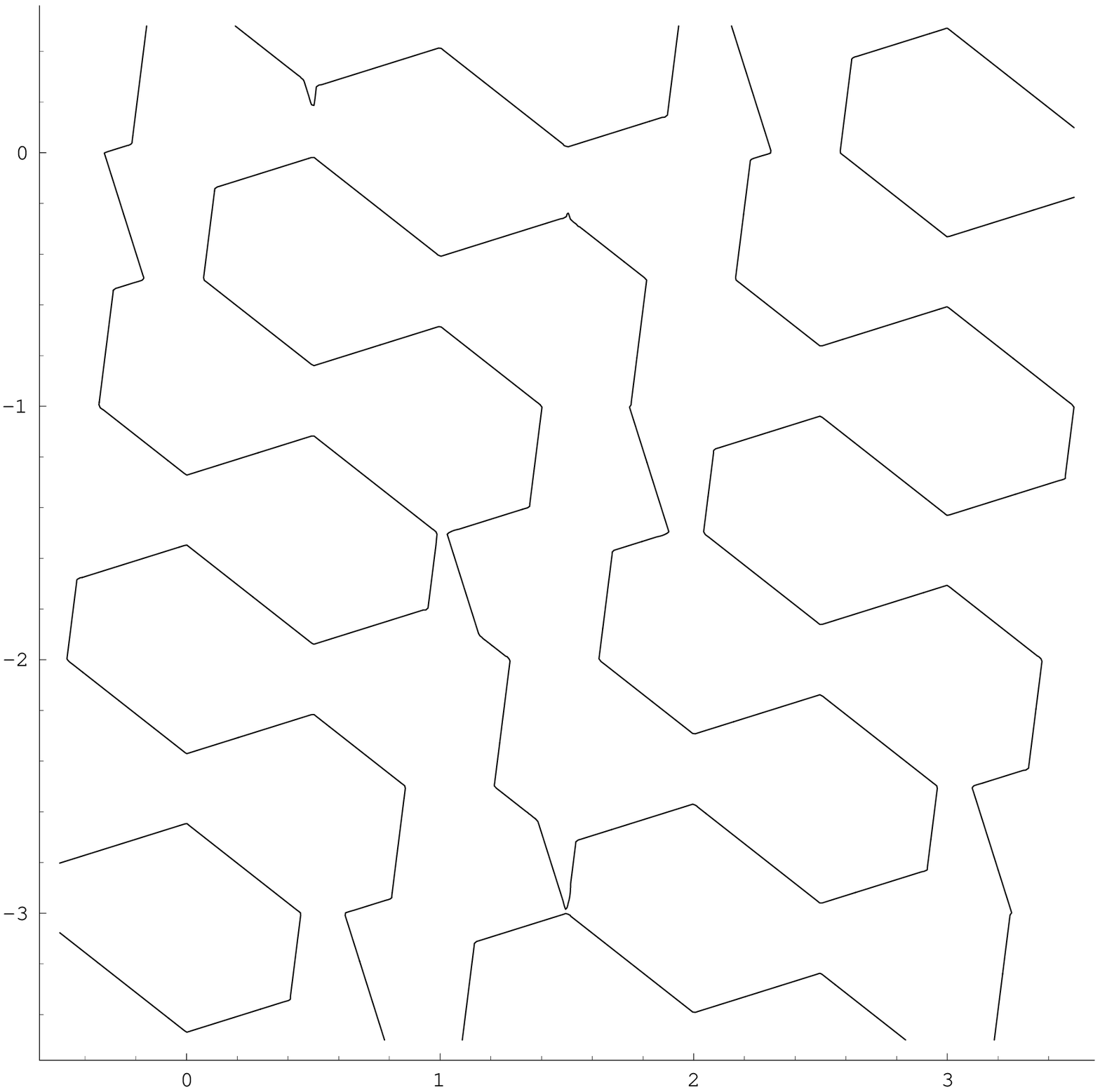}
\cr
\noalign{\vskip3pt}
{\textbf a.}\ \ $\omega = (396,810,1000)$
& 
{\textbf b.}\ \ $\omega = (41,81,100)$
& 
{\textbf c.}\ \ $\omega = (43,82,100)$
\cr 
\noalign{\vskip3pt}
\includegraphics[width=4cm]{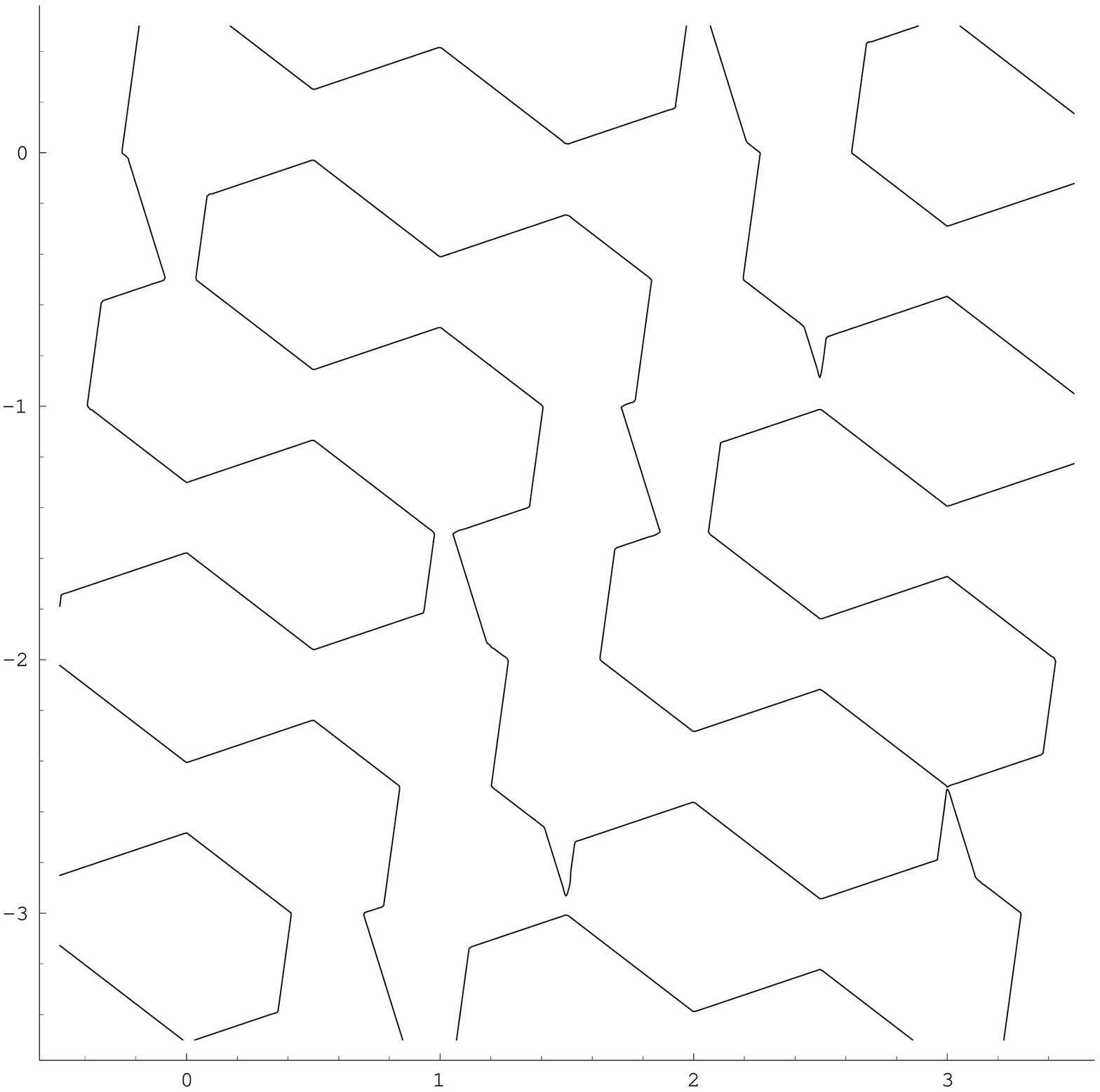}
&
\includegraphics[width=4cm]{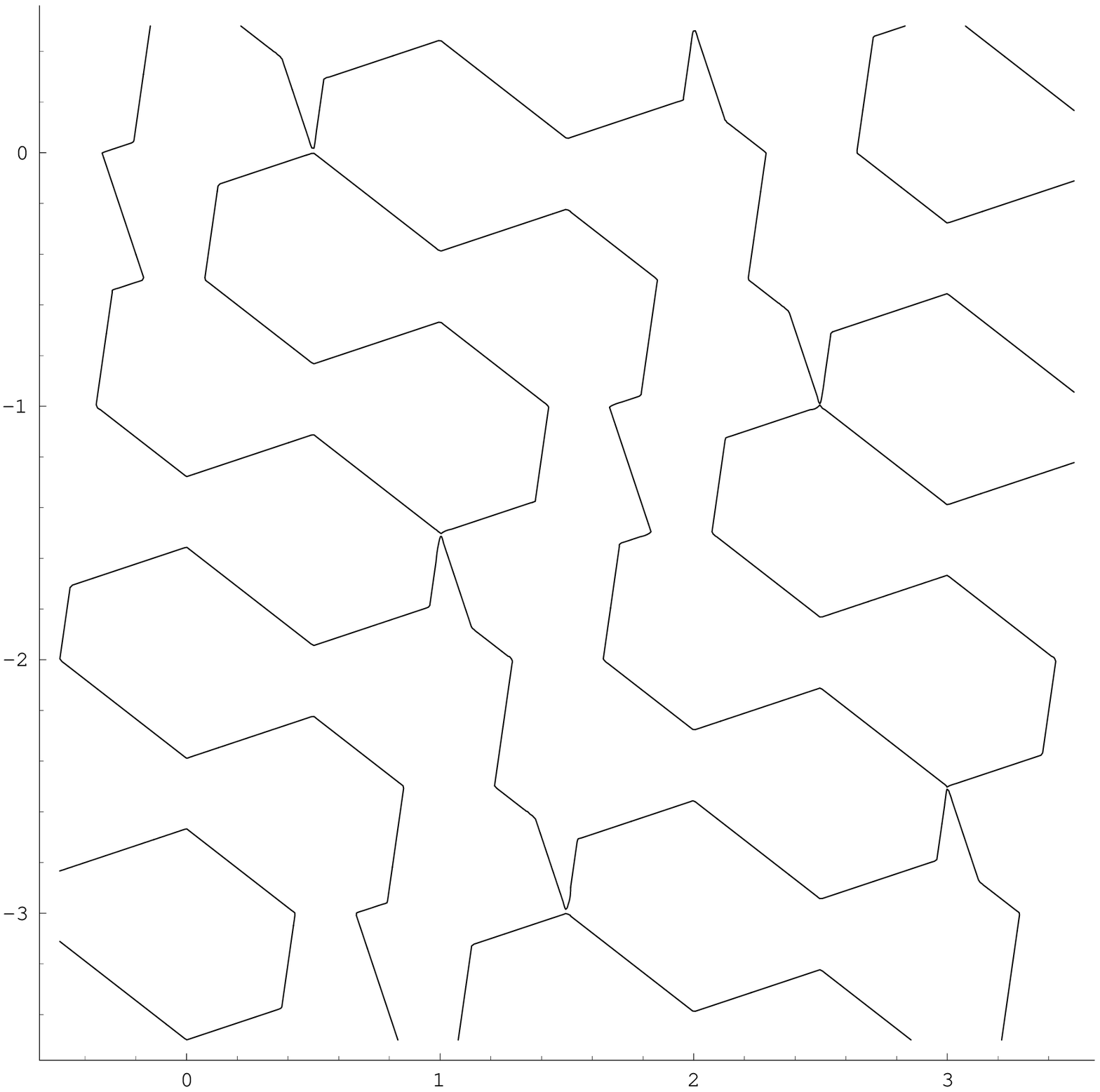}
&
\includegraphics[width=4cm]{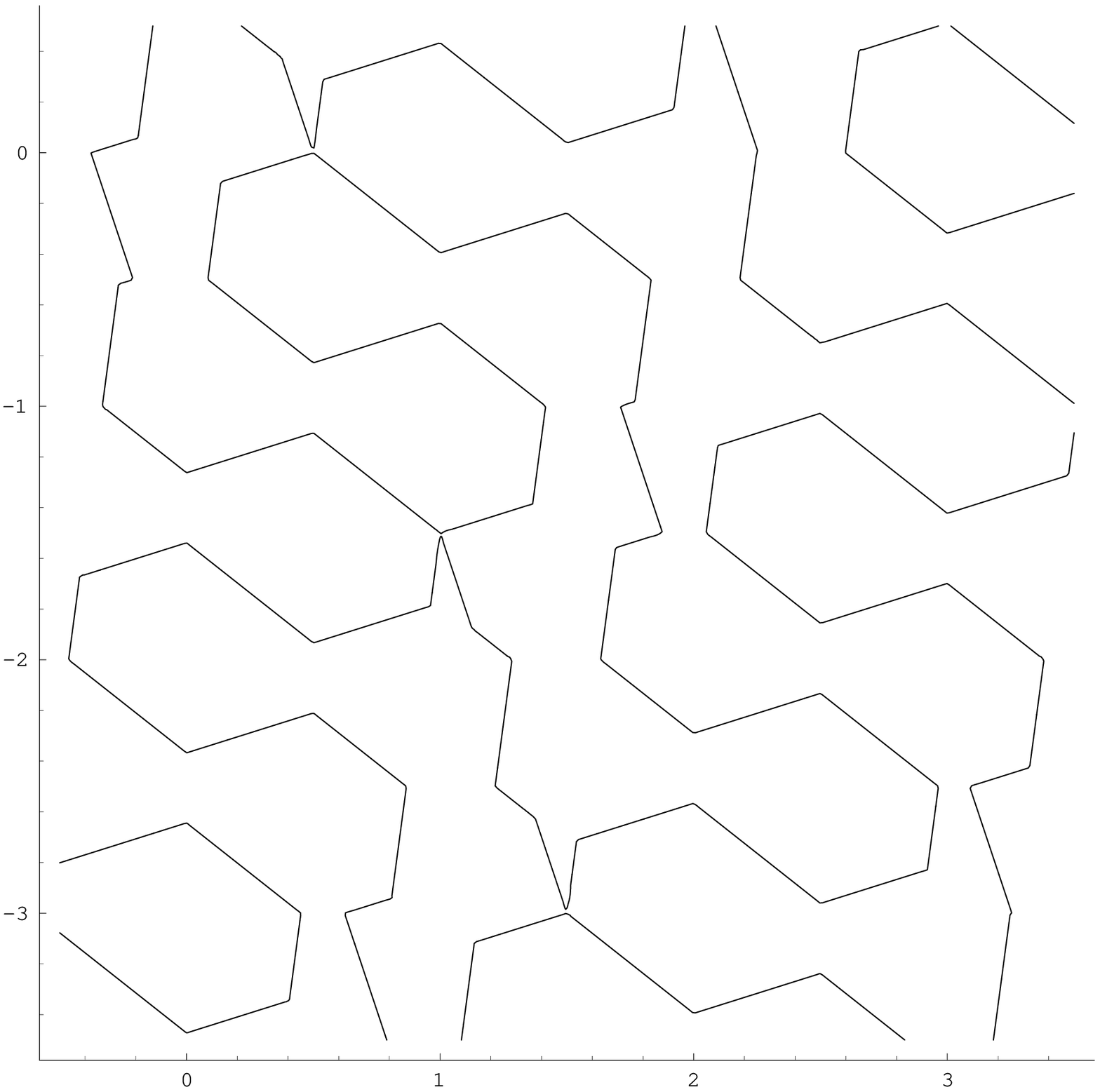}
\cr
\noalign{\vskip3pt}
{\textbf d.}\ \ $\omega = (39,81,100)$
& 
{\textbf e.}\ \ $\omega = (40,80,100)$
& 
{\textbf f.}\ \ $\omega = (43,81,100)$
\cr 
\noalign{\vskip3pt}
\includegraphics[width=4cm]{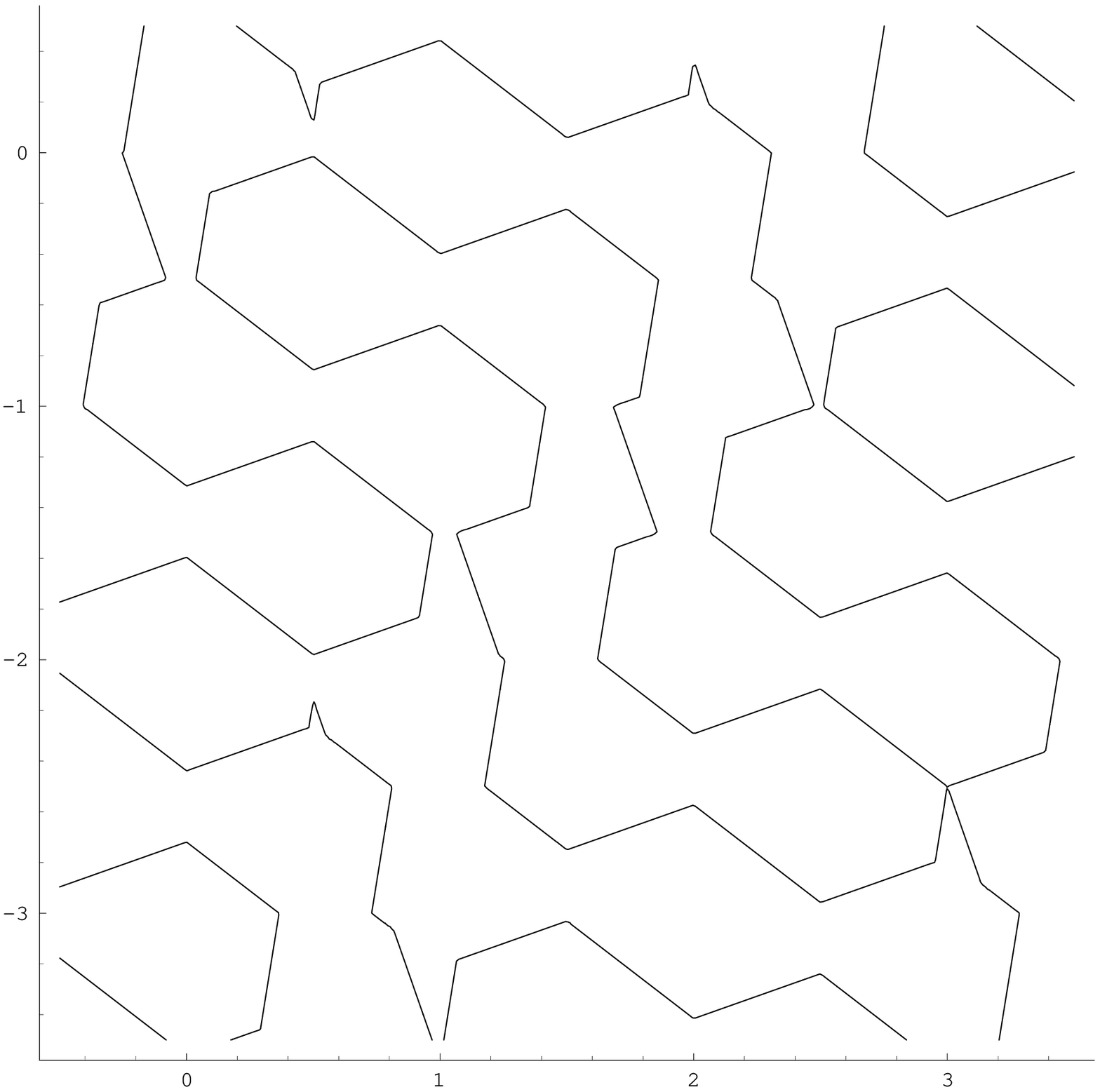}
&
\includegraphics[width=4cm]{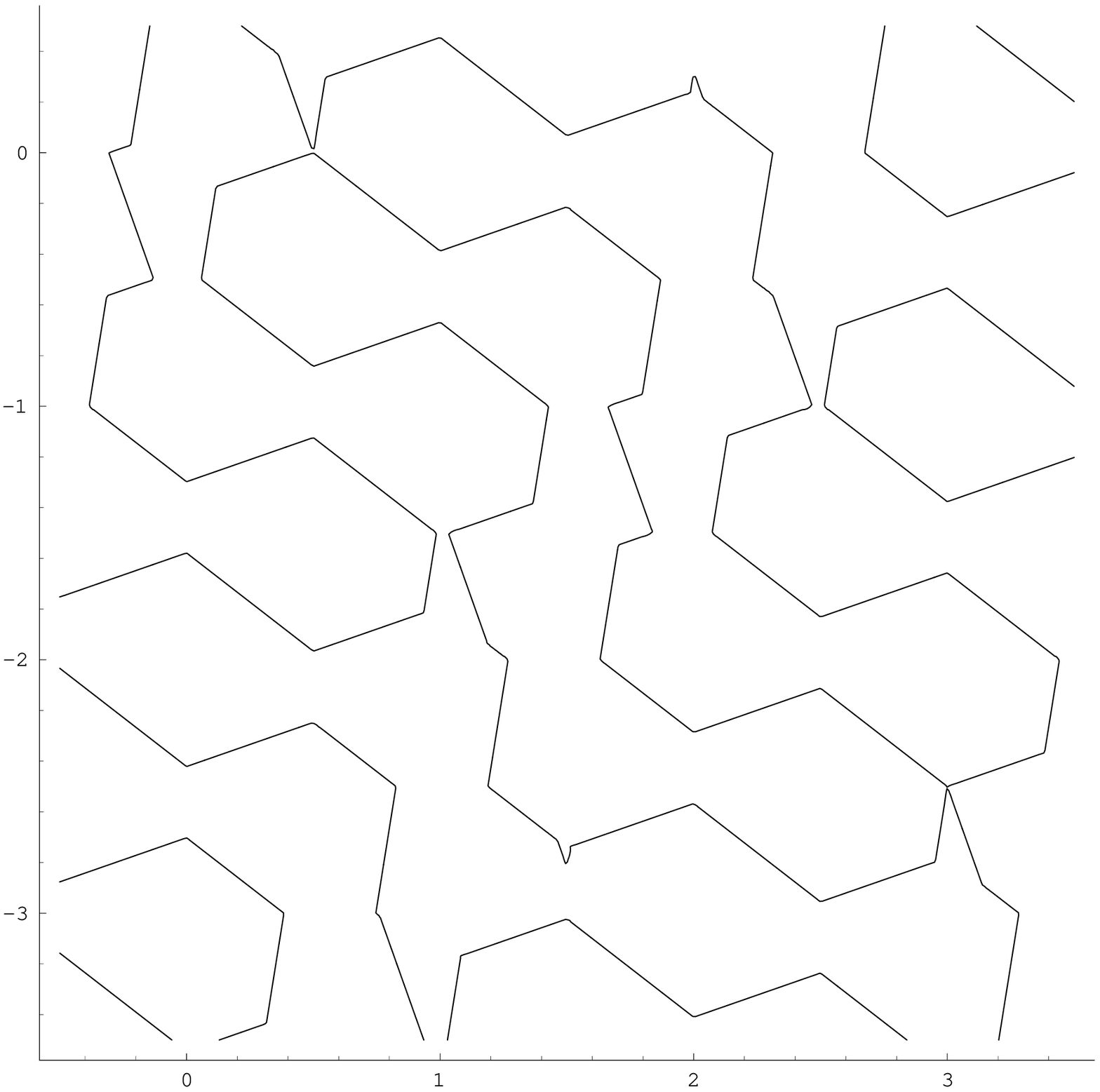}
&
\includegraphics[width=4cm]{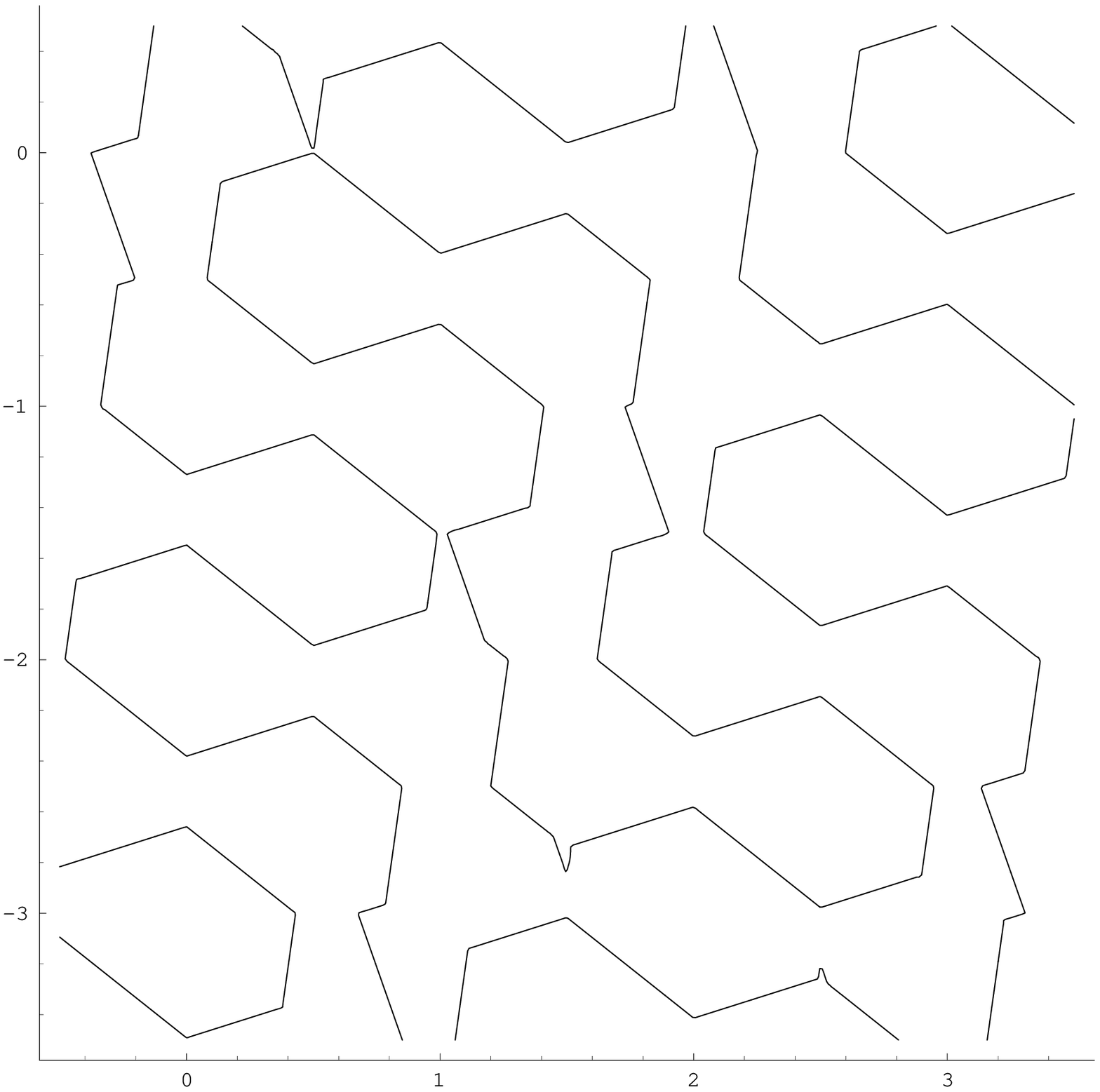}
\cr
\noalign{\vskip3pt}
{\textbf g.}\ \ $\omega = (37,78,100)$
& 
{\textbf h.}\ \ $\omega = (38,78,100)$
& 
{\textbf i.}\ \ $\omega = (43,80,100)$
\cr 
}}
\end{center}
\caption{%
Significant examples of critical sections of $\cP_0$ for $\omega\in\cD_{(2,4,5)}(\cP_0)$.
At the ``center'' of the island ({\textbf e}) there is a saddle connection between the four
critical points (starting from the highest and turning clockwise) $N=p_1$,
$E=p_2+(3,-3,1)$, $S=p_1+(1,-3,2)$ and $W=p_2+(1,-2,1)$. In the subzones $I$, $II$
and $III$ (Fig.~\ref{fig:2.4.5}) the pairs of critical points at the base of the positive
cylinder are, respectively, $N$ and $p_4+(2,-4,1)$ ({\textbf i}), $S$ and $p_3+(2,-2,0)$
({\textbf b},{\textbf c}) and $E$ and $p_3+(-1,-1,0)$ ({\textbf d}, {\textbf g}).
The separating segments correspond respectively to saddle connections between the pairs of
critical points $S$ and $W$ ({\textbf a}), $N$ and $S$ ({\textbf f}) and $N$ and $E$ ({\textbf h}).
}
\label{fig:sections}
\end{figure}

Since the genus is three, only two cylinders may appear and they will be of opposite sign.
The pairs of critical points at the base of cylinders are locally constant; in the square
under investigation the pairings are $p_1,p_4$ for the positive cylinder and $p_2,p_3$ for
the negative one, so that boundaries are always given by an equation like
$<\bH,p_1-p_4+L>=0$. See Figures \ref{fig:2.4.5} and \ref{fig:sections} for a concrete
example worked out in detail.
\par
Finally, a picture of the whole fractal can be obtained through the natural free action on
$\RPt$ of the tetrahedral group $T_d$, whose order is 24 (Fig.~\ref{fig:sm100}(d)).
\begin{figure}
\begin{center}
\includegraphics[width=16cm,bb=120 280 612 732]{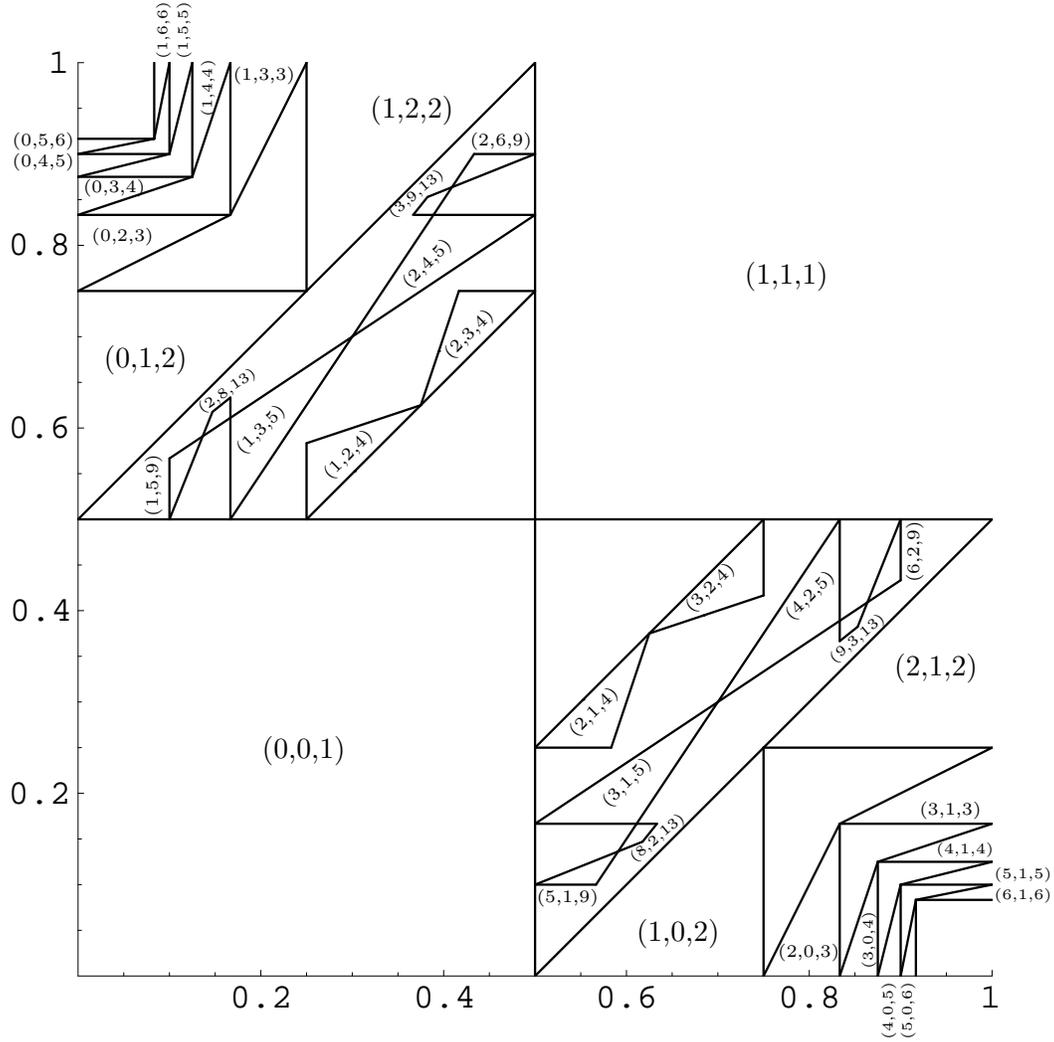}
\end{center}
\caption{%
Analytical boundaries of the biggest islands found in Fig.~\ref{fig:sm100} found using the
algorithm {\textbf A}. All of their boundaries are straight lines segments and the corresponding
label has been reported, when possible, inside the zone itself.
}
\label{fig:sm100an}
\end{figure}
\begin{figure}
\begin{center}
\includegraphics[width=12cm]{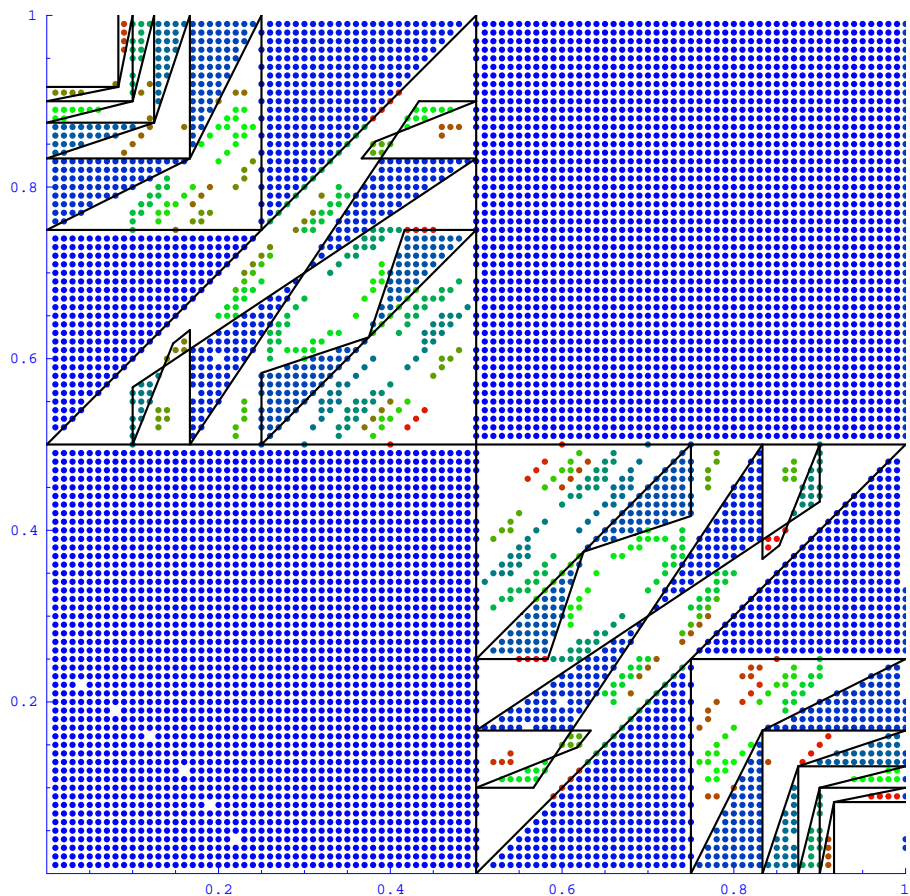}
\end{center}
\caption{%
Comparison between analytical (Fig~\ref{fig:sm100an}) and numerical 
(Fig~\ref{fig:sm100}) boundaries.
}
\end{figure}
\begin{figure}
\begin{center}
\includegraphics[width=12cm]{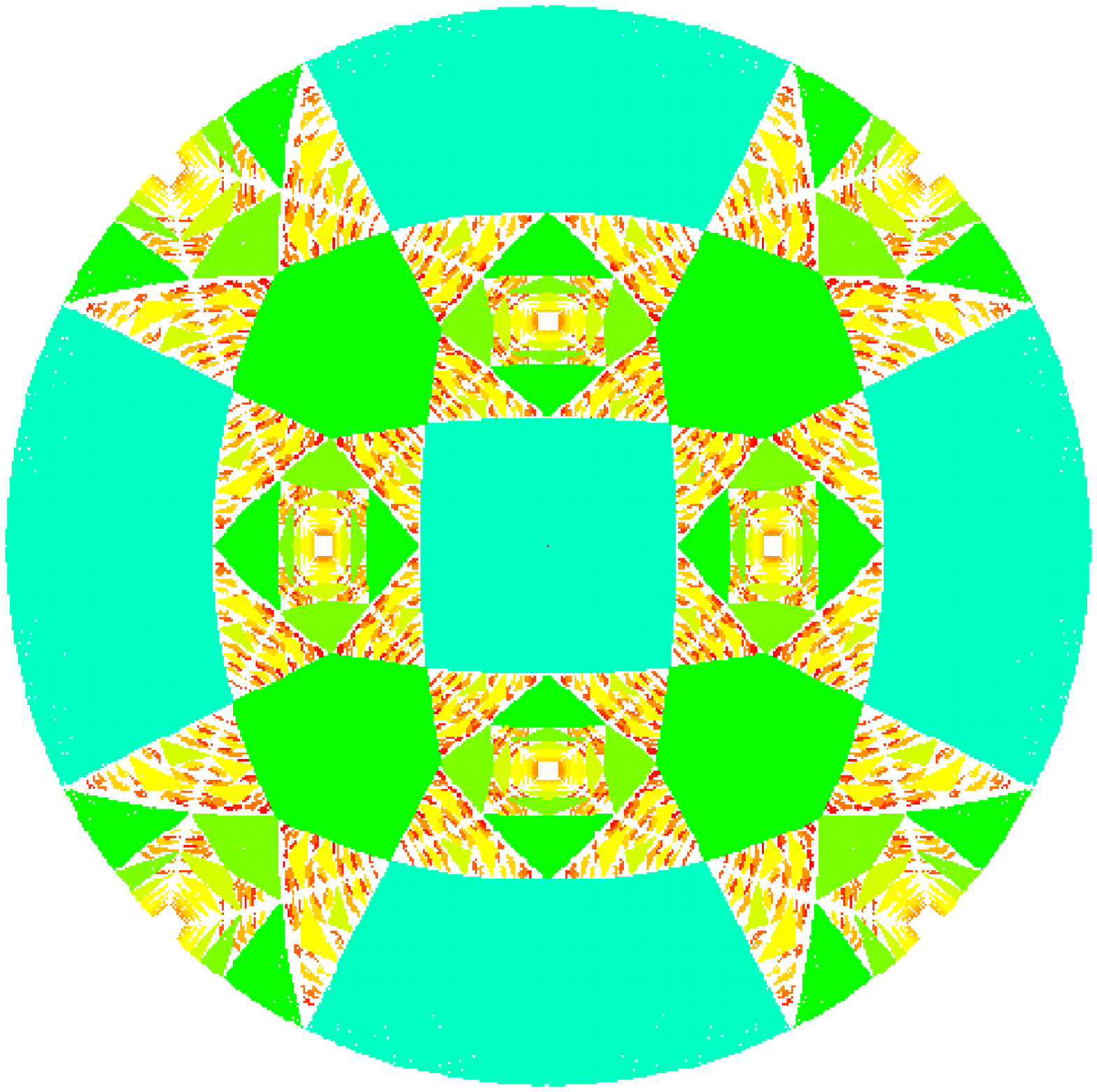}
\end{center}
\caption{%
Fractal image obtained by letting the tetrahedral group $T_d$ act on the square
in Fig.~\ref{fig:sm100} and then projecting on the disc through the stereographic map.
}
\end{figure}
\subsection{Numerical results for $r=10^3$}
\begin{figure}
\includegraphics[width=12cm]{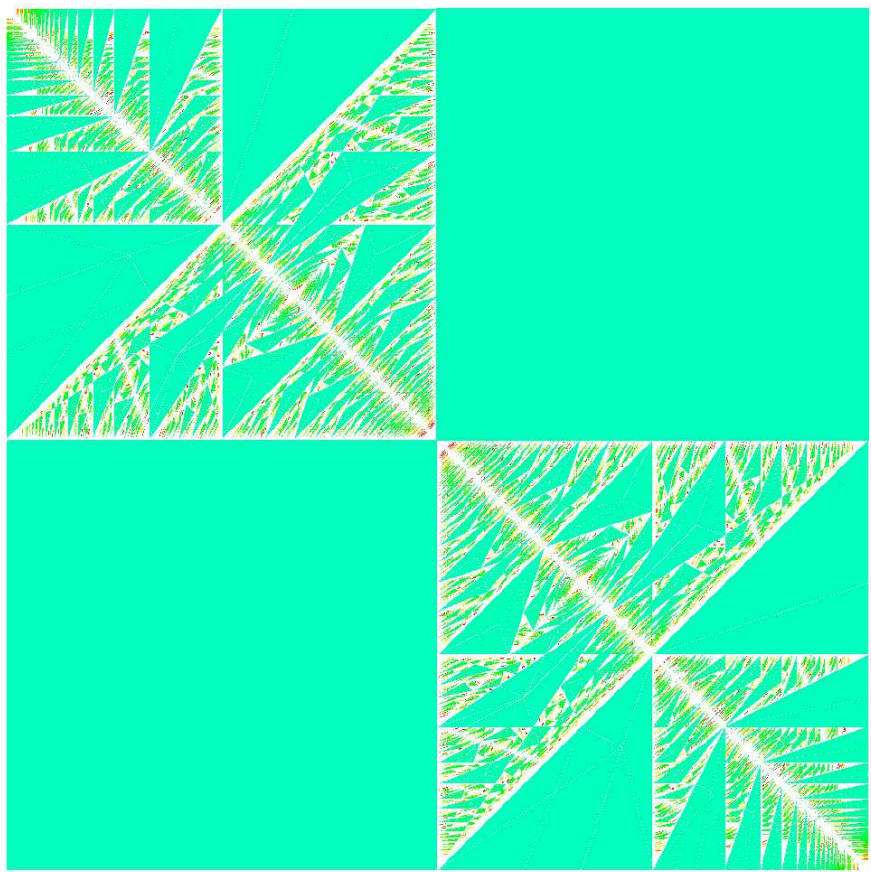}
\caption{%
Numerical plot of the square $[0,1]^2$ of the SM $\cS(\cP_0)=\cS(\mh)$ in the projective
chart $h^z=1$ at a resolution $r=10^3$. The color of the islands goes from green to red as
the norm of the label grows. In the picture are displayed the 1625 islands with at least
four points out of the total 10725 islands found. The missing points that is possible to
see in the interior of some of the islands are due to failure of the numerical algorithm
{\textbf N}, e.g. because of the presence of saddle connections.  The running time for the
$10^6$ steps cycle needed to retrieve these data takes about 1 week on a Pentium 1GHz CPU.
}
\label{fig:sm1000}
\end{figure}

This is the highest resolution reached in~\cite{DL03b}.  In this case for each generic
direction are needed about $10^3$ sections to follow the critical branches, that takes a
time of $~2.5s$ on a $~1GHz$ CPU for the evaluaton of a single label and therefore about
$5\times10^2\times1h\simeq1month$ for sampling the $10^6$ directions of the grid
(Fig.~\ref{fig:sm1000}).
\par
A time of the order of a month is of course rather long, but thanks to the diffusion of
the Linux OS, and therefore of the possibility to build big Unix clusters for cheap using
PCs rather than workstations, this is still not too bad since it is easy to lower by 
a factor 10 the computations time just by dividing the cycle over as many PCs.
This way the running time goes down to just three days, that is a rather acceptable time.
\par
We point out that the situation is radically different in the smooth case: indeed in that
case there is another variable to consider, that is the resolution of the mesh of the
surface, that must be increased together with the grid resolution to avoid errors in the
topology of the curve. The plane sections giving the complete intersection between the
surface and a 2-torus with homology class $(l,m,r)$, $l,m\in[0,r]$, can be as close as
$1/r$ and therefore, if the mesh is too rough, there is the risk that the program will 
jump on the wrong slice. 
\par
Concrete tests show that a mesh resolution of $30$, meaning that the mesh is produced by
dividing the unit cube in a $30^3$ uniform grid, is enough for the $r=10^2$ case but it
must raised to at least $60$ for the $r=10^3$ case, increasing the time for a single label
evaluation to $~15s$, an order of magnitude bigger than in the PL case. This brings back
the time to $~3$ months for the execution, that is indeed the order of the time spent for the
$r=10^3$ calculations made for~\cite{DL03b}.
\par
From the picture~\ref{fig:sm1000} it is rather evident the symmetry with respect to the
anti-diagonal. Apart from this, the picture looks qualitatively very similar to the
pictures found in the previous two cases at the same resolution.
\subsection{Numerical results for $r=10^4$}
\begin{figure}
\begin{center}
\includegraphics[width=4.5cm,bb=200 462 400 722]{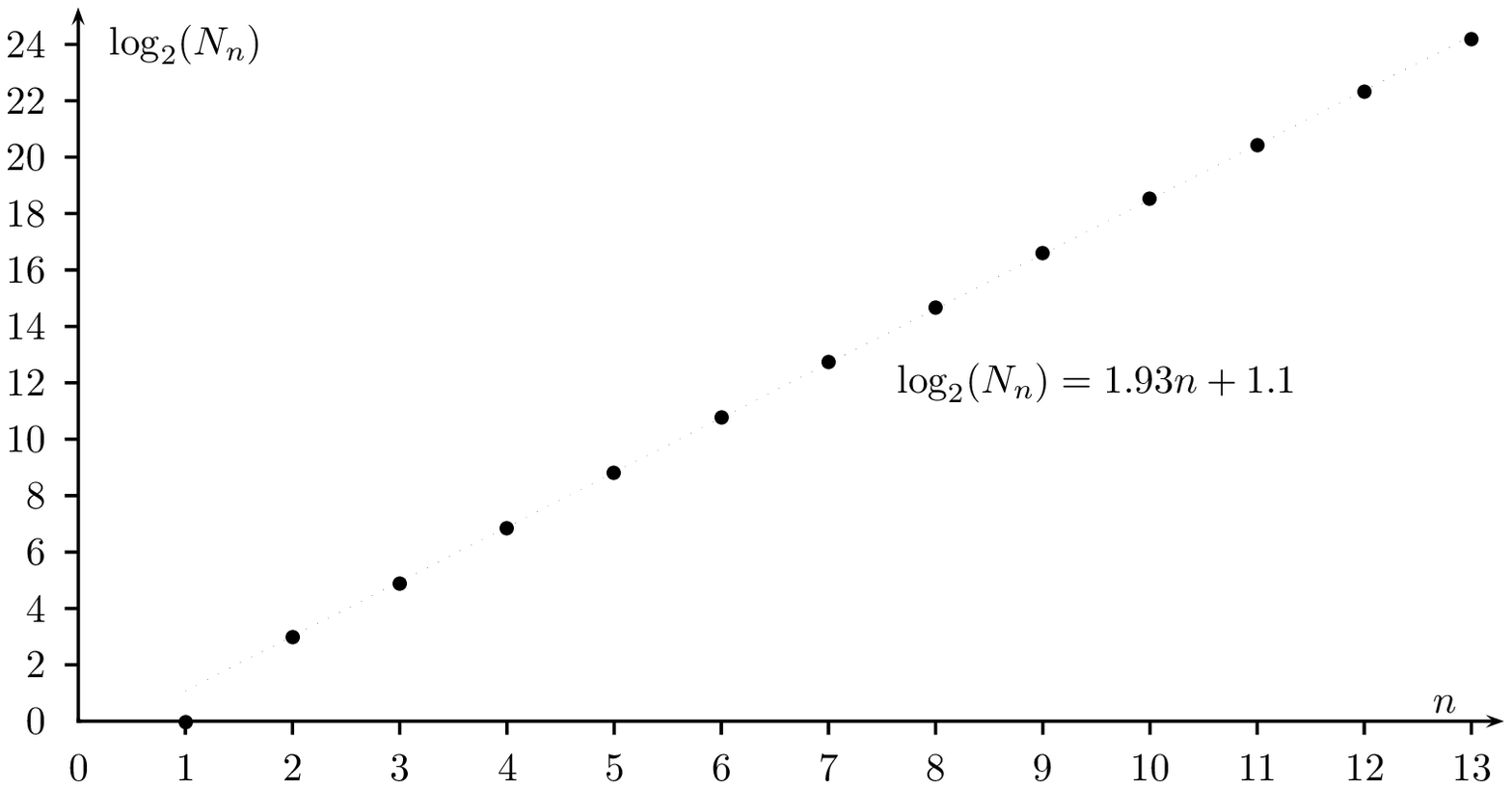}
\end{center}
\caption{%
Box counting evaluation of the Hausdorff dimension of $\cS(\cP_0)$ with the $r=10^4$ data.
$N_n$ is the number of squares of side $2^n$ needed to cover the complement of
$\cS(\cP_0)$, the angular coefficient of the linear fit provides the dimension estimate.
The same evaluation made with the $r=10^3$ data gives a very close result, and also
restricting the fit by canceling a few points at the extremes does not change considerably
the estimate.
}
\label{fig:boxcount}
\end{figure}
Increasing by another order of magnitude the resolution, we increase by an order of magnitude
the number of sections needed to follow a generic leaf, resulting in another factor 5 in the
running time for a single evaluation of a 1-form label, that is now $~10s$, so that the 
total running time on a $1GHz$ CPU reaches $5\cdot10^2\cdot20d=10^4d\simeq30years$.
\par
Such a big running time is rather scary and suggests that there is no hope to go up by
an order of magnitude in resolution without changing some significant algorithm step.
Nevertheless, this big time can be once again brought down to something reasonable
by running the code in 20 PCs and by restricting the numerical analysis to the upper triangle
of the square $[0,.5]\times[.5,1]$ (that reduces computational time by a further factor 8).
Thanks to all these expedients, the running time goes down by two orders of magnitude,
reaching about 3 months, that is indeed about the time that took to us to collect the
$r=10^4$ data.
\par
Note that for the smooth cases this would not be enough, since we must raise the mesh
resolution to $10^2$ and the time for the single evaluation goes up to $50s$; this,
together with the fact that there's a further factor two due to the lack of symmetry with
respect to the antidiagonal, brings the total running time to $~3$ years, definitely not
realisticly affordable.
\par
In Figures
\ref{fig:polSM10000.x0-1000.y9000-10000.eps}-\ref{fig:polSM10000.x4000-5000.y5000-6000.eps} are
shown the numerical results, from which it is clear beyond any doubt that the SM has a
fractal-like self-repeating structure, even though no explicit construction is known for
it.
\par
\begin{figure}
\begin{center}
\includegraphics[width=8.5cm,bb=100 530 400 722]{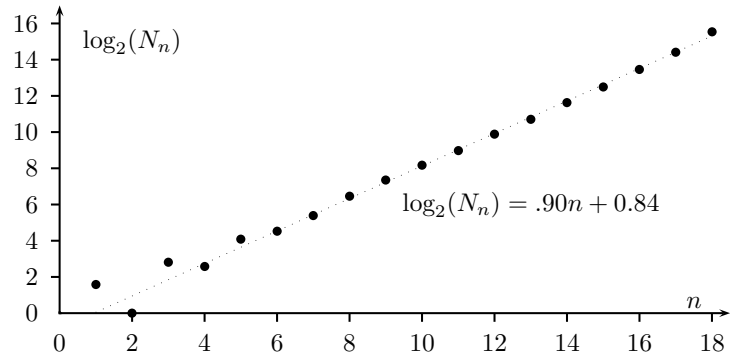}
\end{center}
\caption{%
``Area distribution'' evaluation of the Hausdorff dimension of $\cS(\cP_0)$ with the
$r=10^4$ data. $N_n$ is the number of islands whose area is between $2^{-n}$ and
$2^{-n-1}$. The estimate of the Hausdorff dimension is given by the double of the angular
coefficient of the linear fit and varies by $\pm0.2$ by restricting even by little the
number of points on which the fit is made.
}
\label{fig:areasdistr}
\end{figure}

Numerical evaluations of the Hausdorff dimension $d_{\cP_0}$ of the fractal set, namely
the complement of $\cS(\cP_0)$ in $\RPt$, has been performed using the box counting
technique (Fig.~\ref{fig:boxcount}) and a sort of ``area distribution'' technique
(Fig.~\ref{fig:areasdistr}).
\par
The first, and more reliable, technique involves partitioning the square $[0,1]^2$ in
$2^n$ identical squares and evaluating the number of squares needed to cover the fractal.
With the $r=10^3$ data, $n$ can get as big as $\log_2(10^3)\simeq10$ and a linear fit
gives an evaluation of $d_{\cP_0}\simeq1.93$. The $r=10^4$ data allow $n$ to go up to
$13$, representing the deeper results on the Hausdorff dimension of such fractals to date;
figure~\ref{fig:boxcount} shows that the scaling law is linear to a high degree of
accuracy and a linear fit gives again an estimate of $d_{\cP_0}\simeq1.93$, making us
rather confident in the accuracy of this numerical result.
\par
\begin{figure}
\begin{center}
\includegraphics[width=5cm,bb=200 362 400 622]{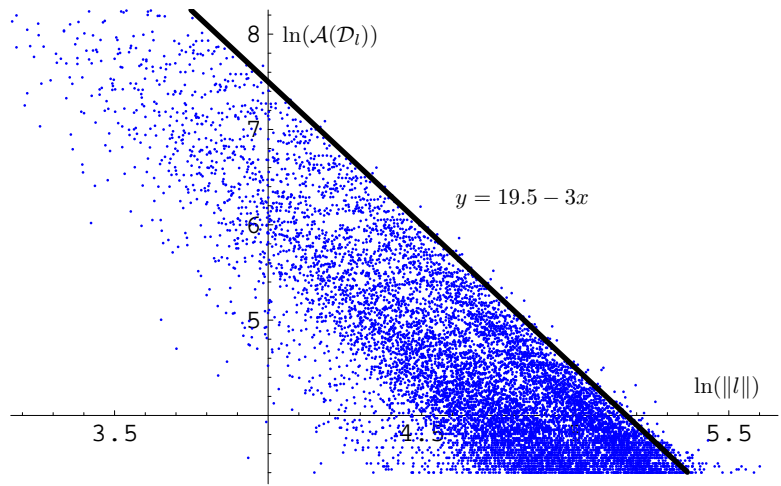}
\end{center}
\caption{%
Logarithmic plot of the islands' areas versus the norm of the corresponding label.
There is a very good agreement between the numerical data and the Dynnikov conjecture
claiming that $\cA(\cD_l)\leq C/\|l\|^3$.
}
\label{fig:dynConj}
\end{figure}
The second one involves counting the number of zones whose size is between $b^n$ and
$b^{n+1}$, where $b>1$. Computations with sevaral small values, between 2 and 1, were
performed and all of them gives rough estimate of $d_{\cP_0}\simeq1.76$, further than
expected from the evaluation given by the box counting. It is not clear to us the reason
of this disagreement, that did not manifest for the smooth cases~\cite{DL03a}, but it is
not impossible that it may be due simply to the fact that this evaluation stabilizes at
higher resolutions and therefore this evaluation is at the current state of things more or
less unreliable (as simple numerical tests testify).
\par
Finally, the $r=10^4$ data provide enough detail to test numerically a conjecture by
Dynnikov~\cite{Dyn99} and give a graphical evidence of theorem~\ref{thm:labels}.
The Dynnikov conjecture claims that the area of the islands satisfy a relation
$\cA(\cD_l(\cM))\leq C/\|l\|^3$ for some positive real number $C$ depending only on the
polyhedron$\cM$; the numerical data suggests that the exponent 3 cannot be improved any further  
(see Fig.~\ref{fig:dynConj}).
\par
\begin{figure}
\includegraphics[width=12cm]{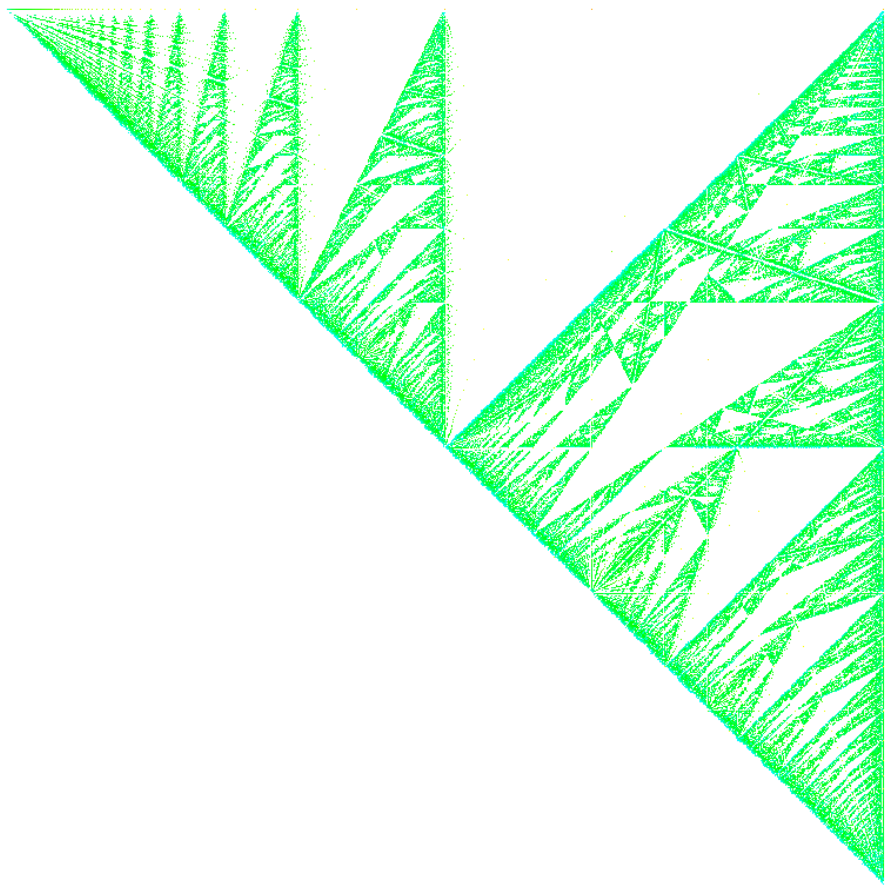}
\caption{%
Plot of the labels of the 494041 islands found in the numerical analysis, at the
resolution $r=10^4$, of the portion of the SM $\cS(\cP_0)$ contained in the triangle of
vertices $(0,1)$, $(.5,1)$ and $(.5,.5)$ in the projective chart $h^z=1$. According
to theorem~\ref{thm:labels}, the closure of the set of labels is equal to the complement
of the interiors of the islands; the striking closeness of the two pictures is one of the
best indirect tests of the correctness of the implementation of the {\textbf N} algorithm in
our C++ library NTC.
}
\label{fig:sm1000labels}
\label{fig:labels}
\end{figure}
According to theorem~\ref{thm:labels}, the closure of the set of labels of $\cS(\cP_0)$, 
seen as points of $\RPt$, is equal to the complement of the interior of the islands.
At this resolution about $5\cdot10^5$ different labels are found and their image
in $\RPt$ (Fig.~\ref{fig:labels}) is one of the best indirect checks of correctness
of the library NTC.
\section{Conclusions}
We proved in this paper that the structure of foliations induced on triply periodic
embedded polyhedra by constant 1-forms is identical to the one induced on smooth
surfaces, and we believe that the same line of arguments can also be used to further
generalize the theorems to embedded piecewise smooth surfaces.
This fact extends to ``Morse PL functions'' the association of a SM; in the most
interesting cases such SM have a fractal nature (this condition is true for an open subset
of all triply periodic PL functions, e.g for all triply periodic functions close enough to
$\mh$).
\par
Surfaces satisfying property SP1 (see section~\ref{sec:sp}) are particularly rich, their
SM being equal to the SM of any function having them as level set. We exploited this fact
by studying the case of the triply periodic surface $\cP_0$ obtained extending in the
three coordinate directions a truncated octahedron. The simplicity of the triangulation of
the surface allowed us to improve by an order of magnitude, respect to the results
obtained in~\cite{DL03b}, the resolution of the numerical analysis of the fractal
and, as a consequence, to perform several numerical tests on conjectures and theorems.
\section{Acknowledgments}
We want to thank first of all S.P. Novikov for introducing the subject. We are also deeply
in debt with S.P. Novikov and I.A. Dynnikov for their interest in our work and for several
precious and fruitful scientific suggestions and discussions that were essential for the
present work. In particular, we thank I.A. Dynnikov for suggesting the analytical algorithm
to retrieve the islands boundaries.
All numerical calculations were made with computers kindly provided by the INFN section
of Cagliari (www.ca.infn.it) and CRS4 (www.crs4.it). We also acknowledge financial support 
from the Cagliari section of INFN and from the Departments of Physics and of Mathematics
of the University of Cagliari.
\bibliography{pol}
\vfill\eject
\begin{figure}
\includegraphics[width=11cm,bb=200 410 400 722]{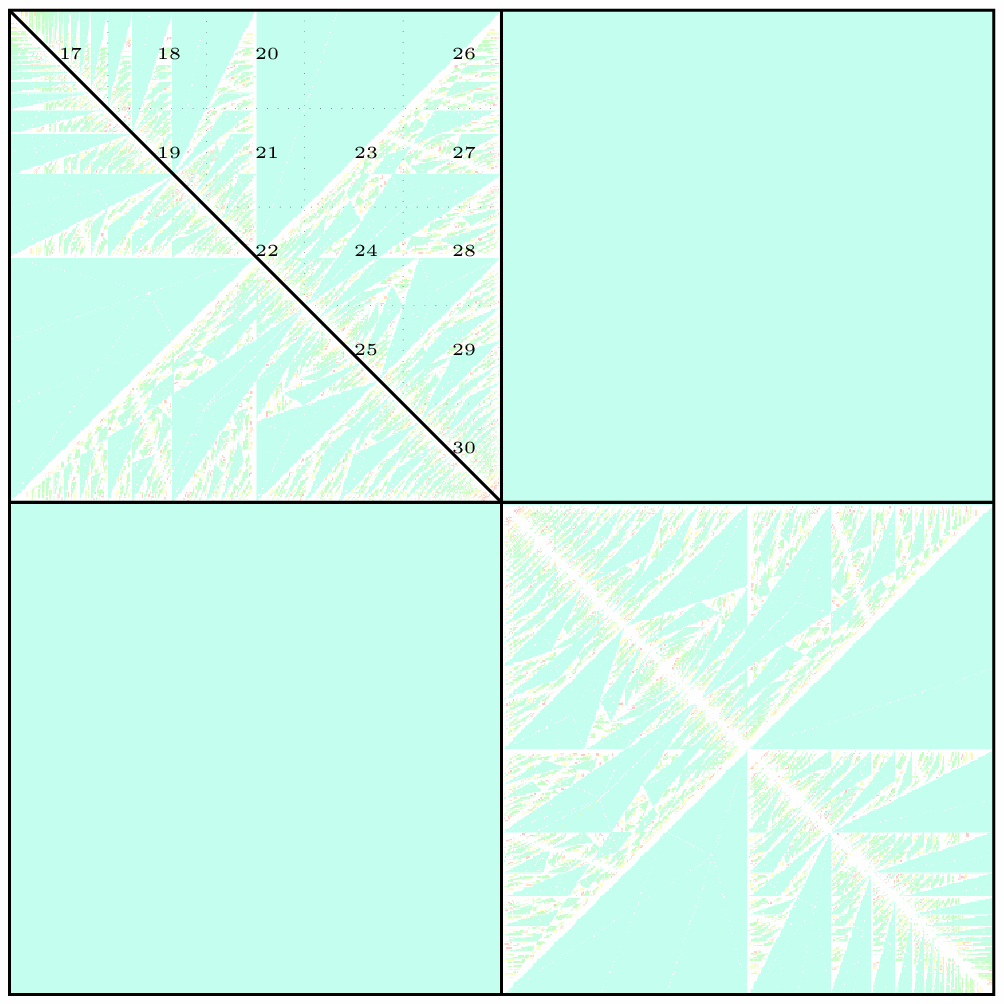}
\caption{%
The next 14 pictures show, in full detail, the fractal structure of $\cS(\cP_0)=\cS(\mh)$
sampled at the resolution $r=10^4$. In order to minimize the CPU time, the numerical
analysis was limited to the triangle $[(0,1,1)]$, $[(1/2,1,1)]$ $[(1/2,1/2,1)]$, divided
above in smaller squares and triangles by a uniform grid and labeled by the corresponding
figure number; the full picture can be retrieved by applying the ``topological'' symmetry
characteristic of this surface, namely the symmetry with respect to the antidiagonal, and
then the 24 symmetries coming from the action of the tetrahedral group $T_d$. To retrieve these 
data we used about 20 1GHz Pentium III CPUs for about 3 months.In all
pictures, the color of the islands goes from green to red as the norm of the label
grows. Note that the square $[.4,.5]\times[.9,1]$ is fully contained inside the island
$\cD_{(1,2,2)}$ and therefore its picture will not be shown.
}
\label{fig:sm10000}
\end{figure}
\vfill\eject
\fig{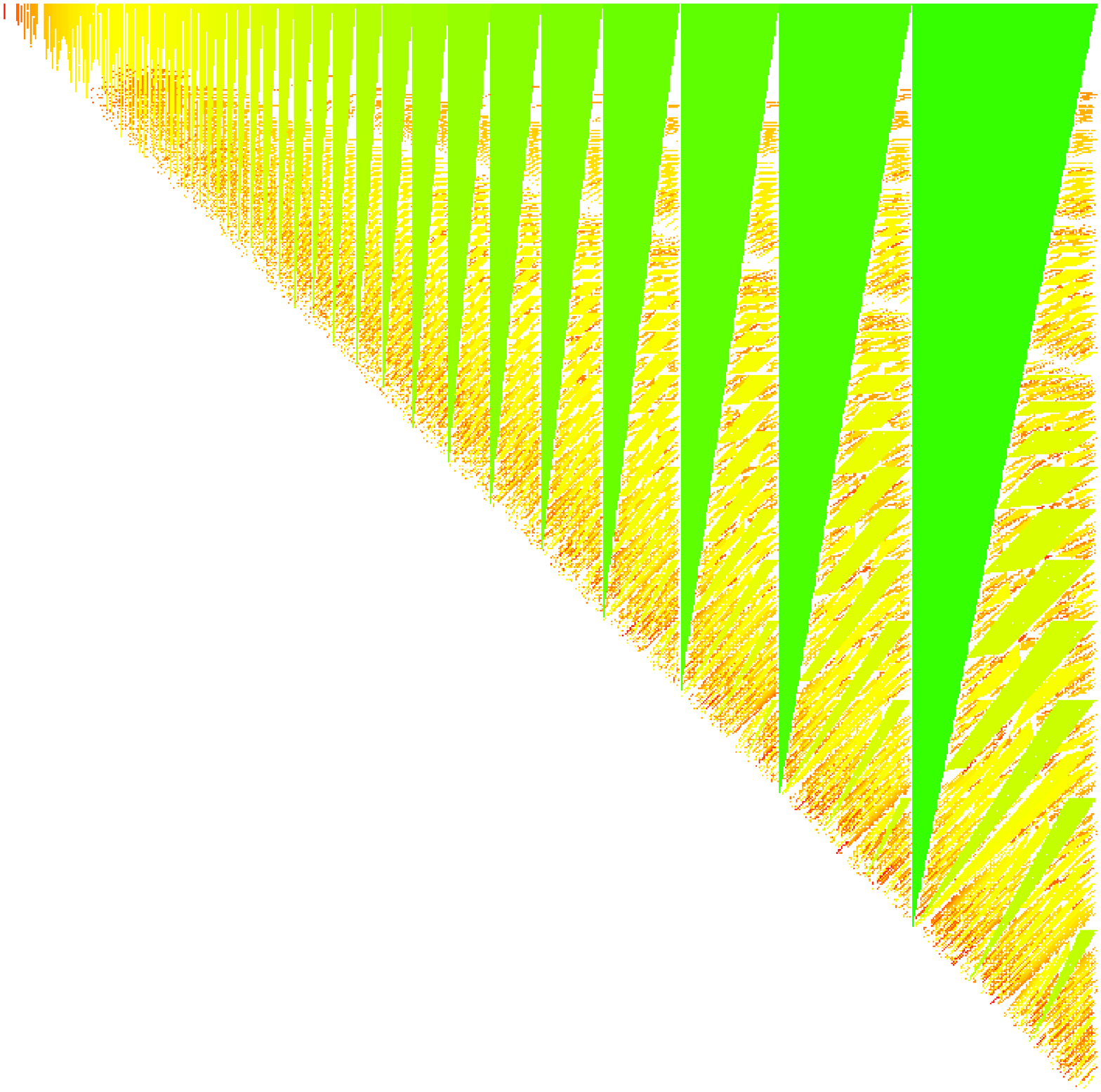}{Detail of $\cS(\cP_0)$ in $[0,.1]\times[.9,1]$}
\vfill\eject
\fig{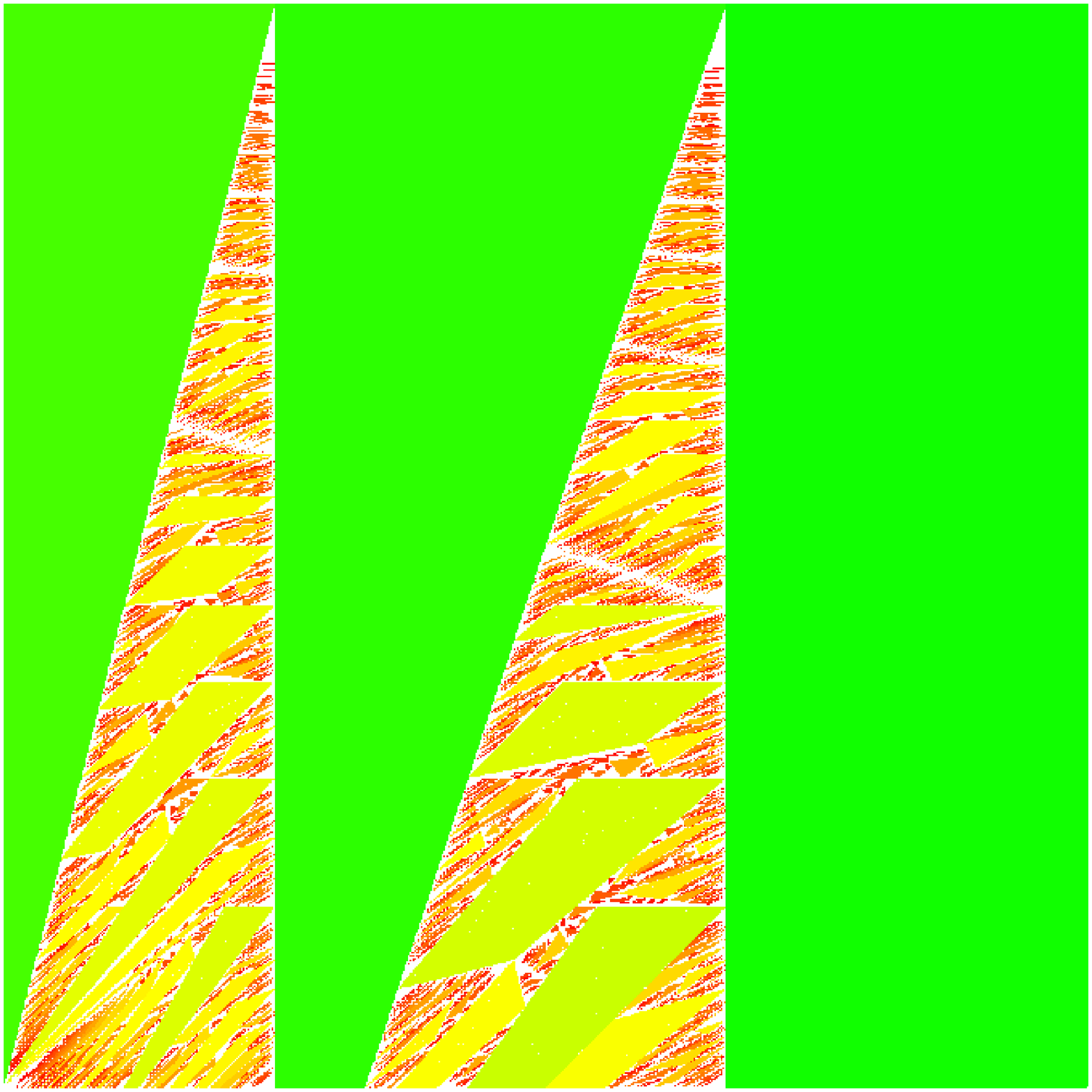}{Detail of $\cS(\cP_0)$ in $[.1,.2]\times[.9,1]$}
\vfill\eject
\fig{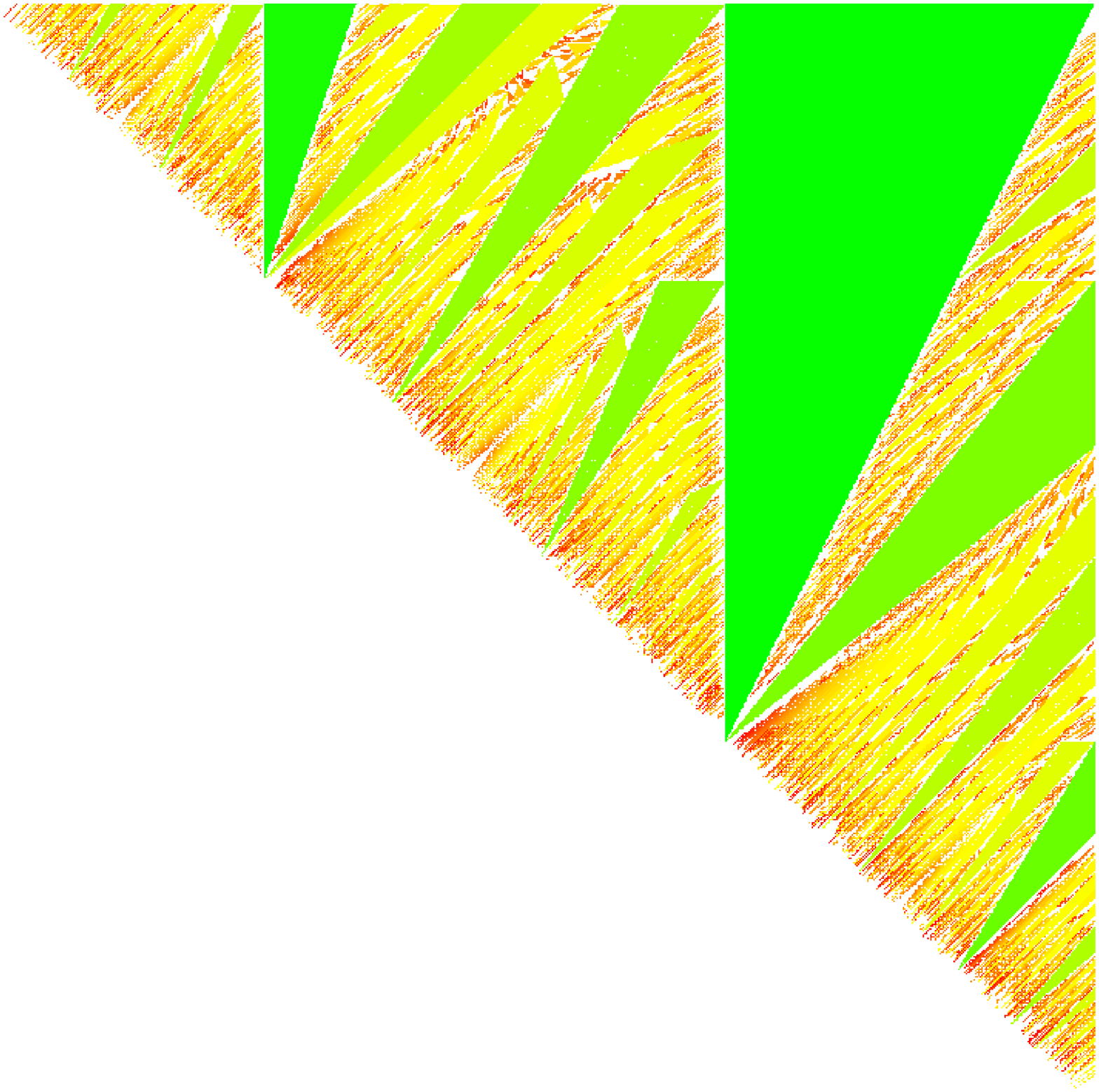}{Detail of $\cS(\cP_0)$ in $[.1,.2]\times[.8,.9]$}
\vfill\eject
\fig{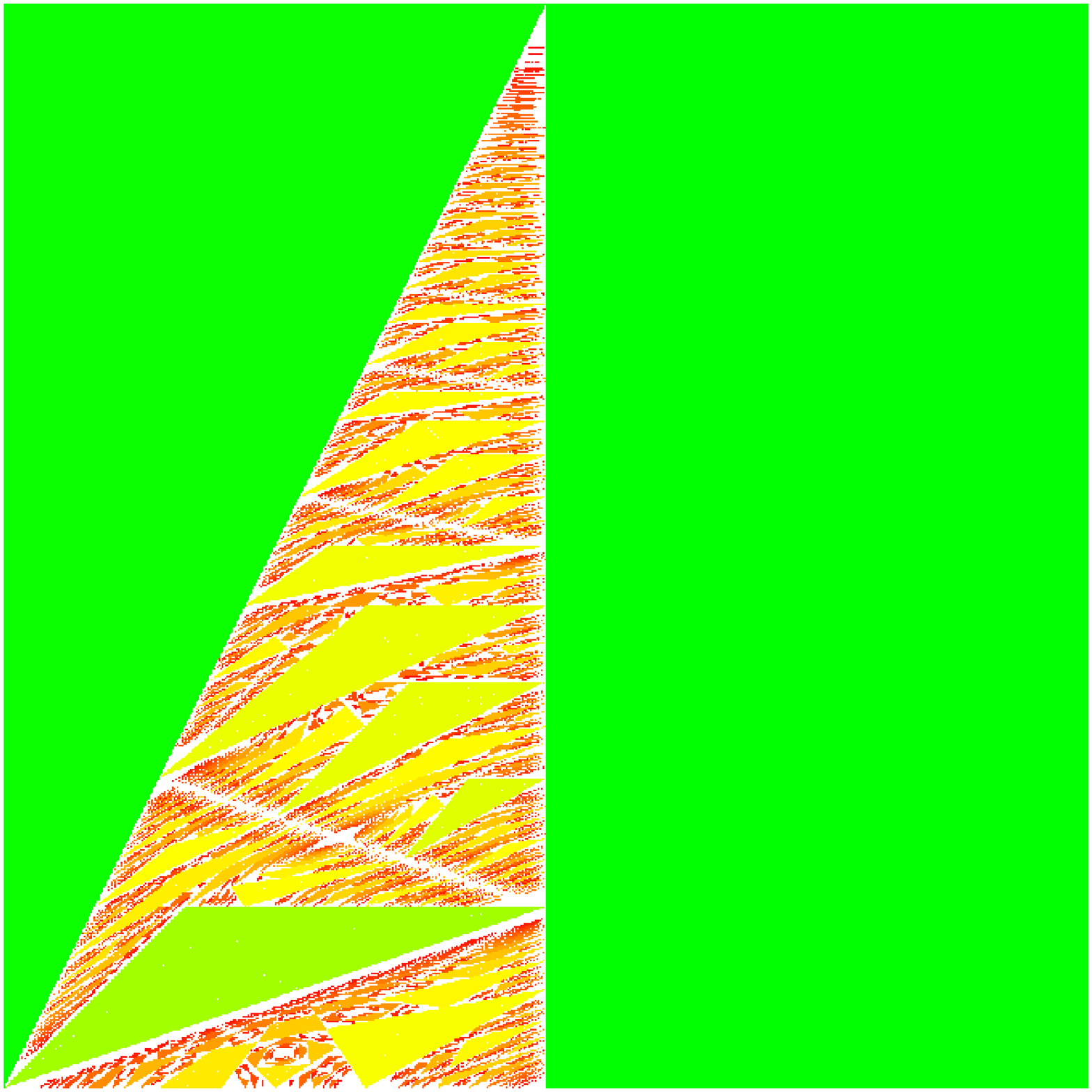}{Detail of $\cS(\cP_0)$ in $[.2,.3]\times[.9,1]$}
\vfill\eject
\fig{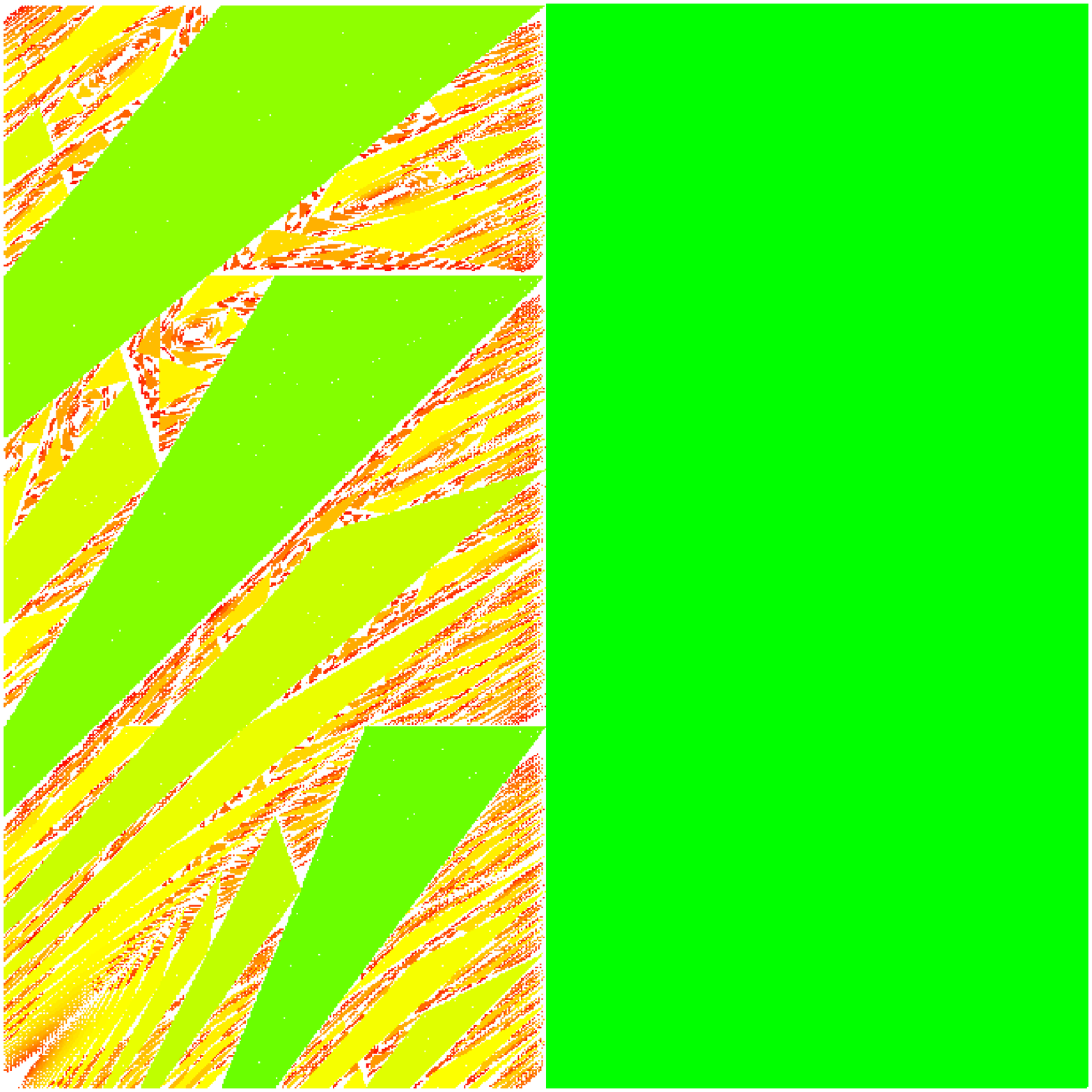}{Detail of $\cS(\cP_0)$ in $[.2,.3]\times[.8,.9]$}
\vfill\eject
\fig{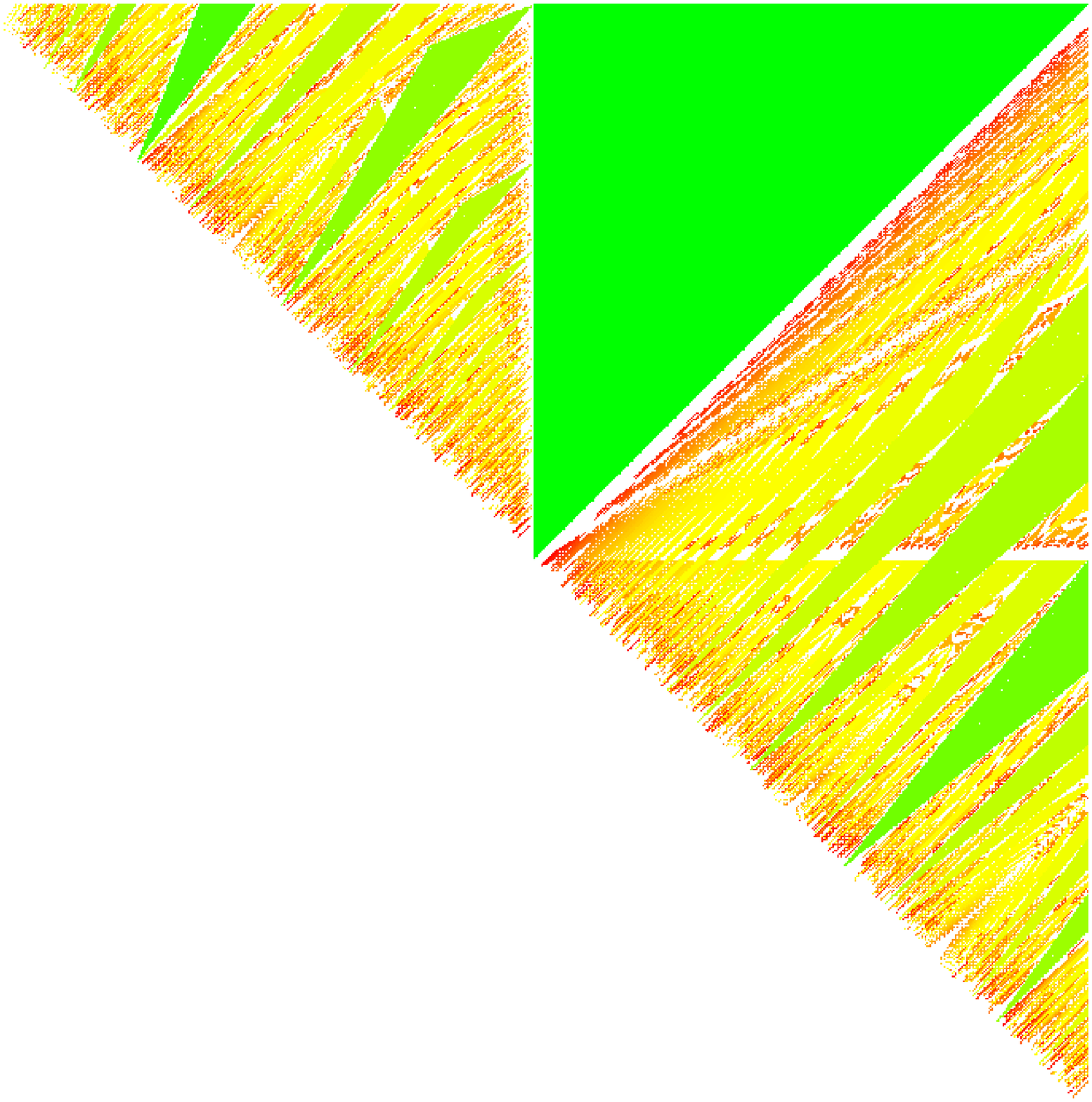}{Detail of $\cS(\cP_0)$ in $[.2,.3]\times[.7,.8]$}
\vfill\eject
\fig{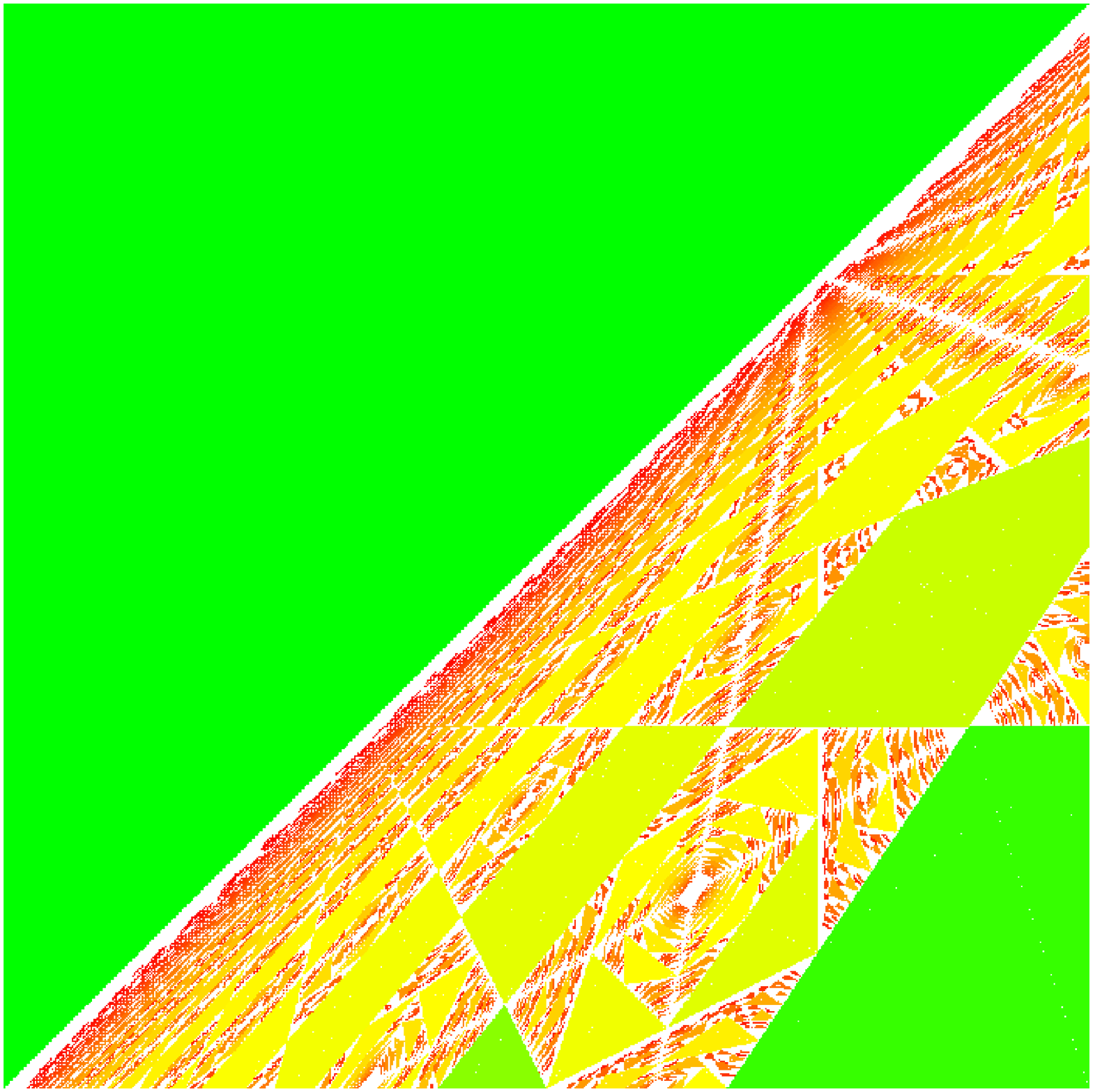}{Detail of $\cS(\cP_0)$ in $[.3,.4]\times[.8,.9]$}
\vfill\eject
\fig{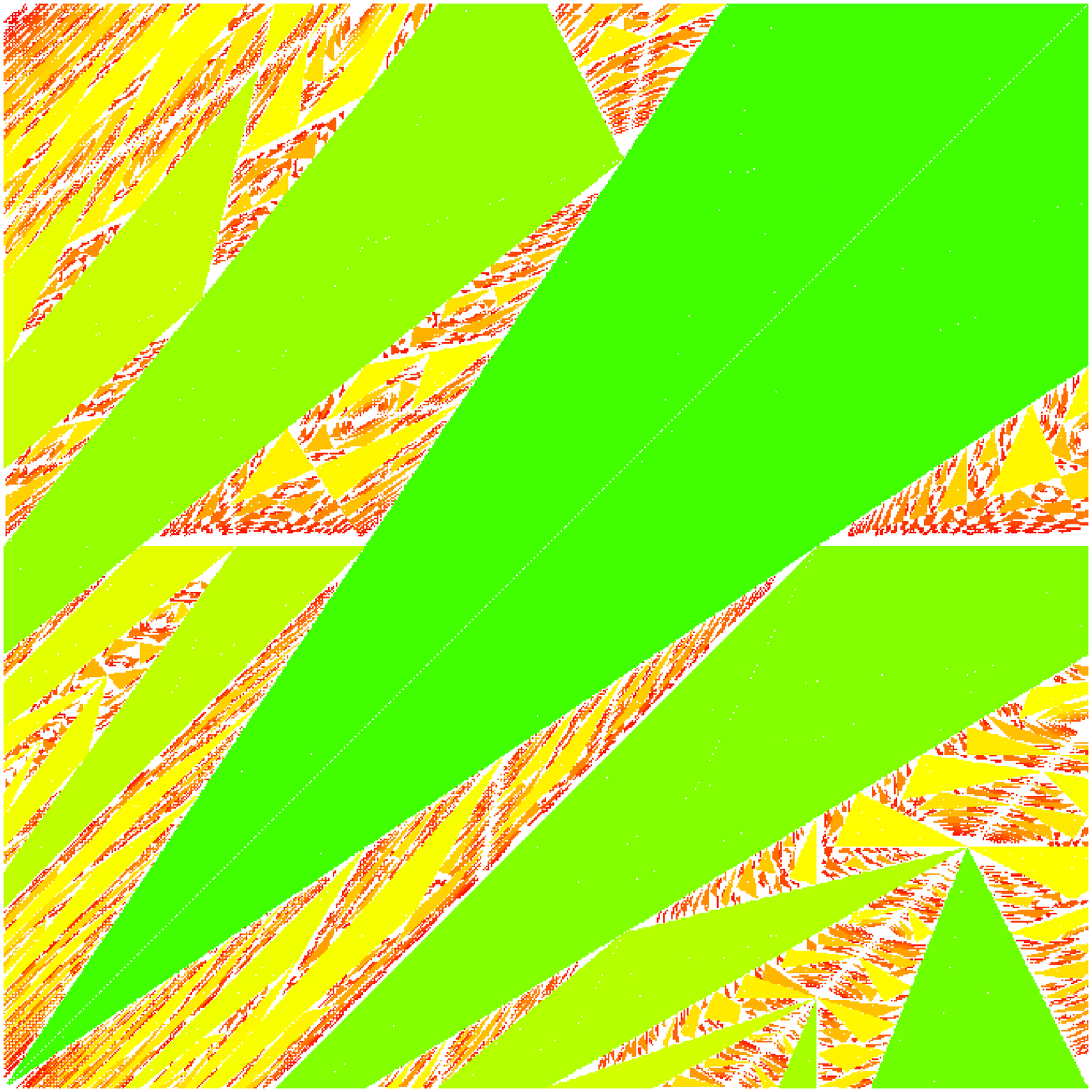}{Detail of $\cS(\cP_0)$ in $[.3,.4]\times[.7,.8]$}
\vfill\eject
\fig{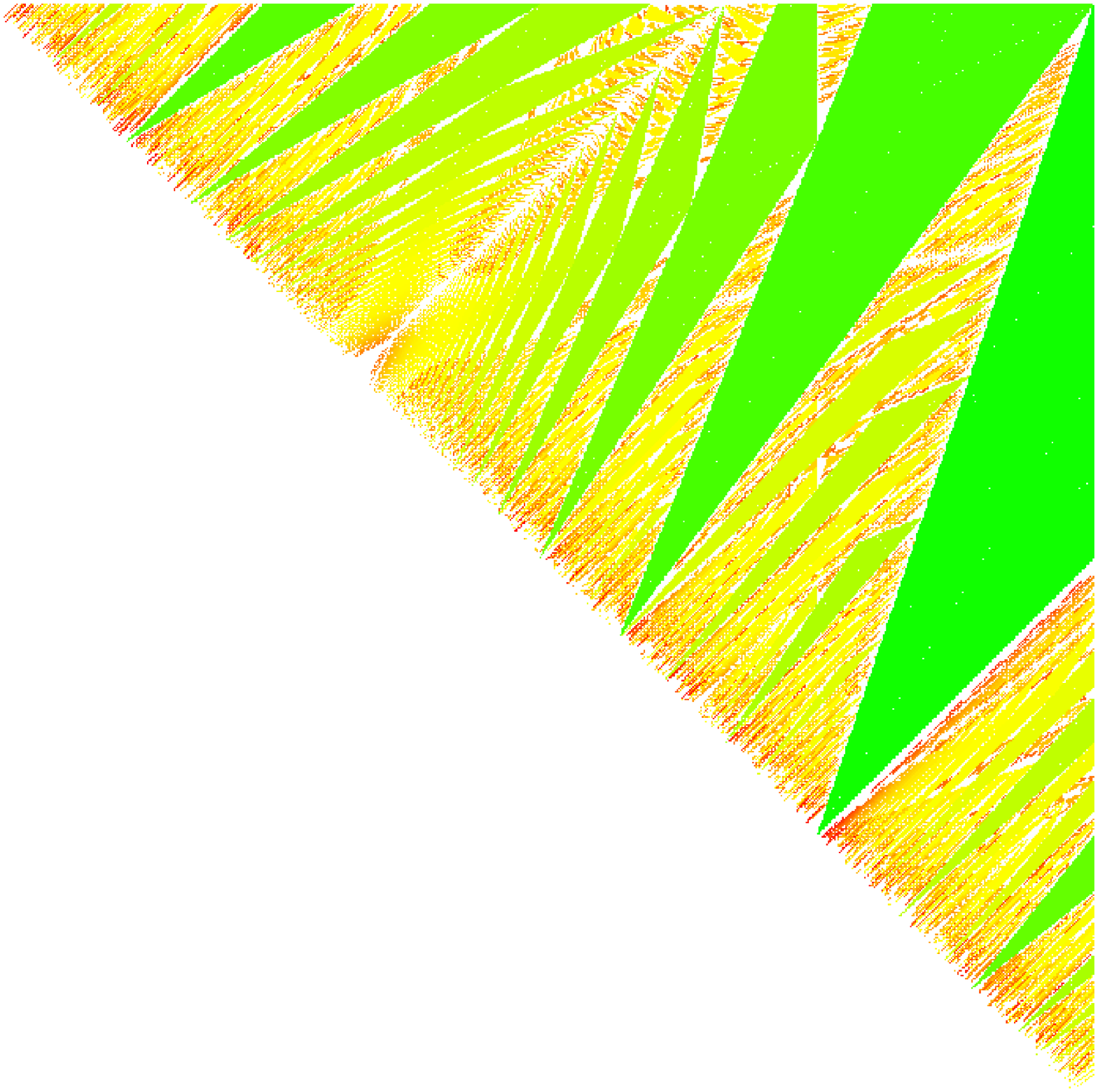}{Detail of $\cS(\cP_0)$ in $[.3,.4]\times[.6,.7]$}
\vfill\eject
\fig{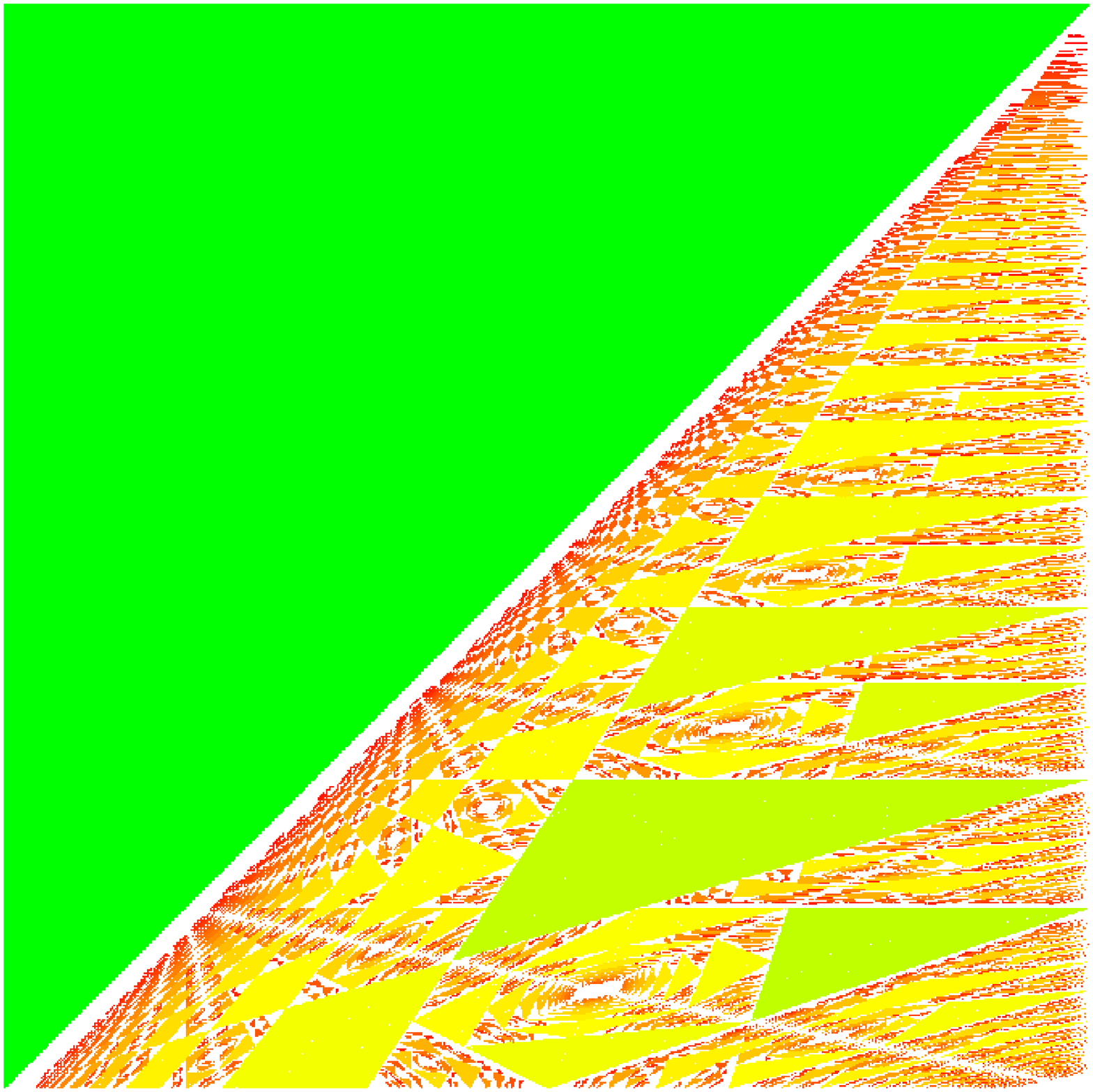}{Detail of $\cS(\cP_0)$ in $[.4,.5]\times[.9,1]$}
\vfill\eject
\fig{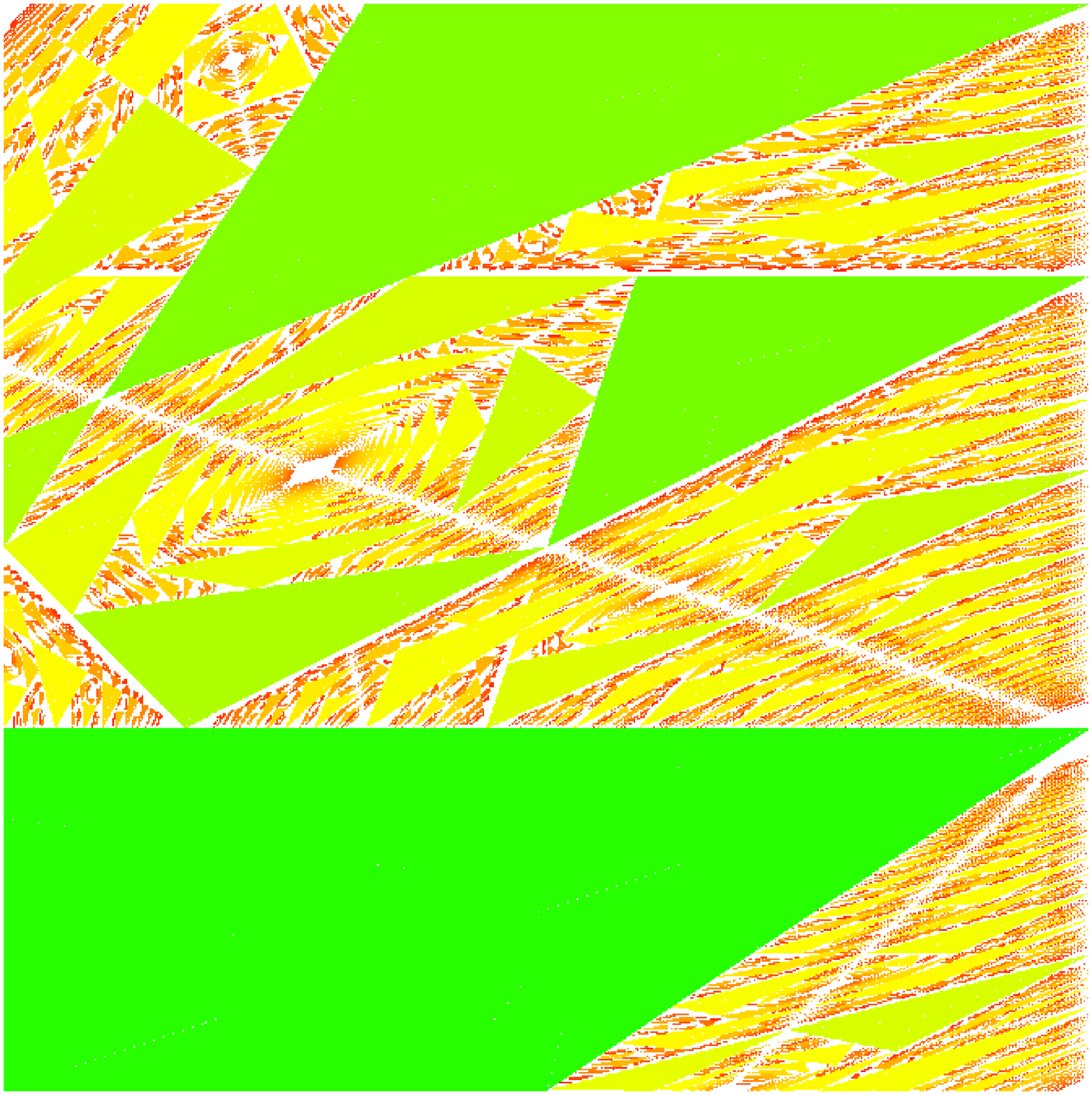}{Detail of $\cS(\cP_0)$ in $[.4,.5]\times[.8,.9]$}
\vfill\eject
\fig{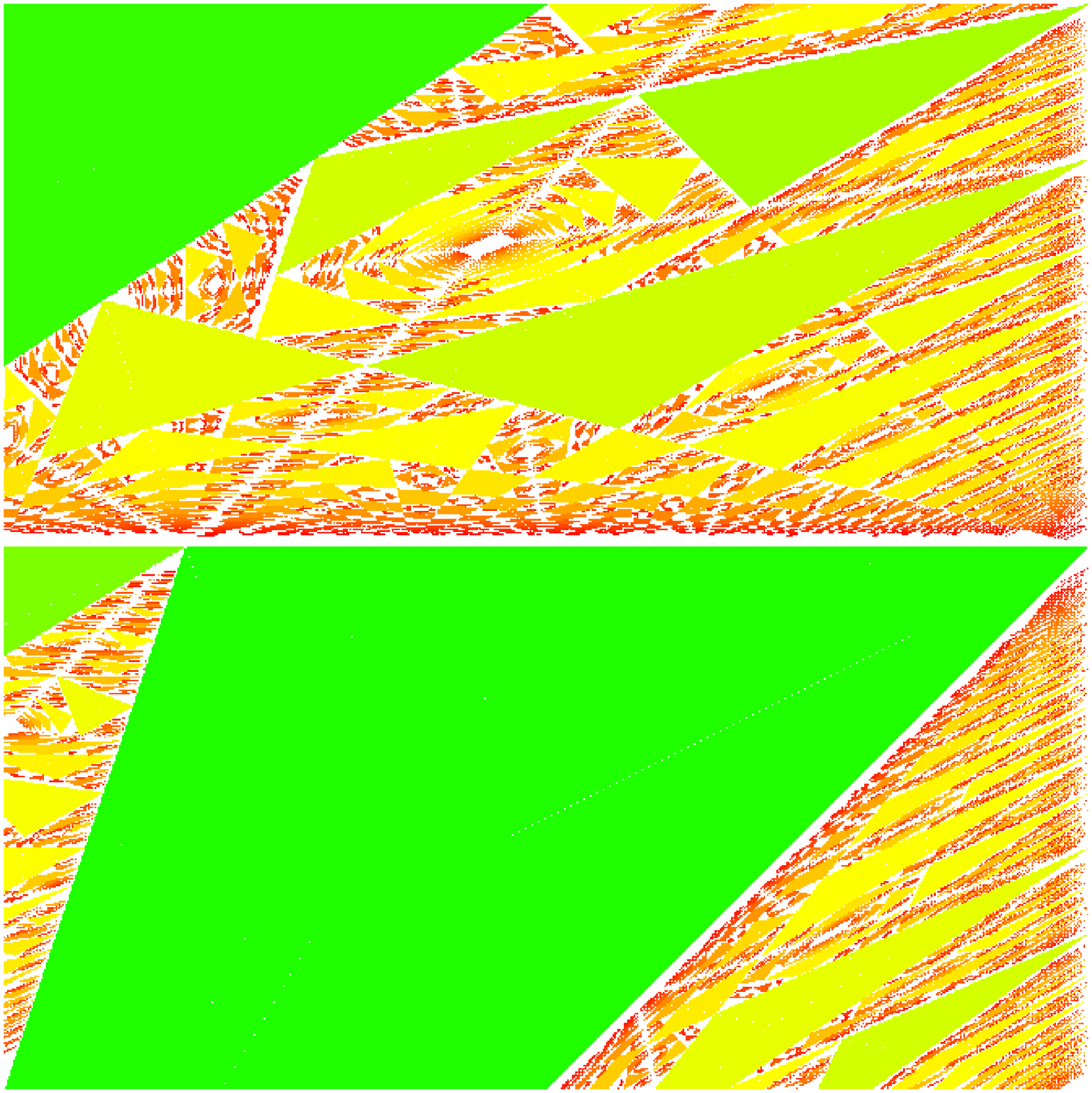}{Detail of $\cS(\cP_0)$ in $[.4,.5]\times[.7,.8]$}
\vfill\eject
\fig{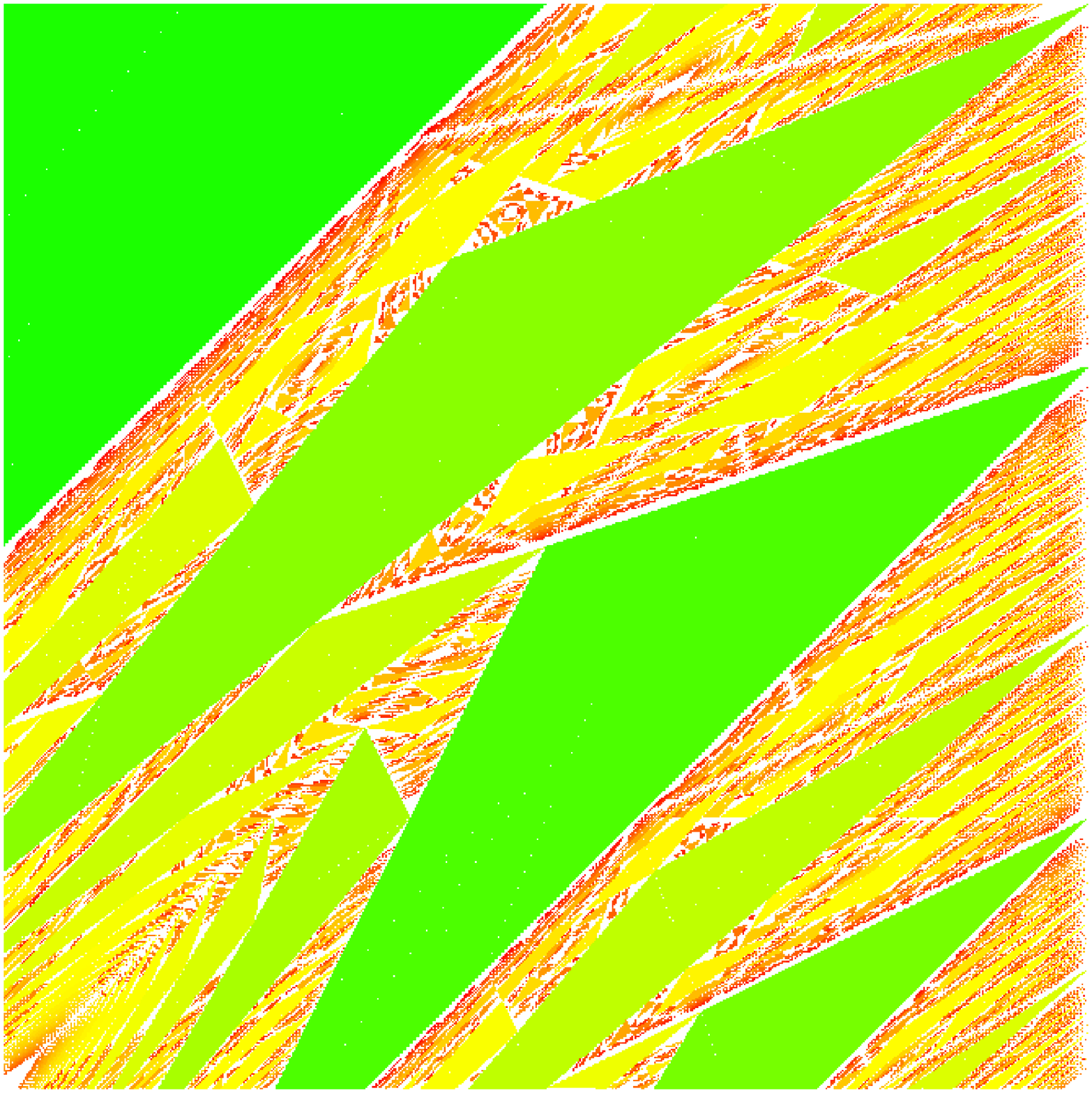}{Detail of $\cS(\cP_0)$ in $[.4,.5]\times[.6,.7]$}
\vfill\eject
\fig{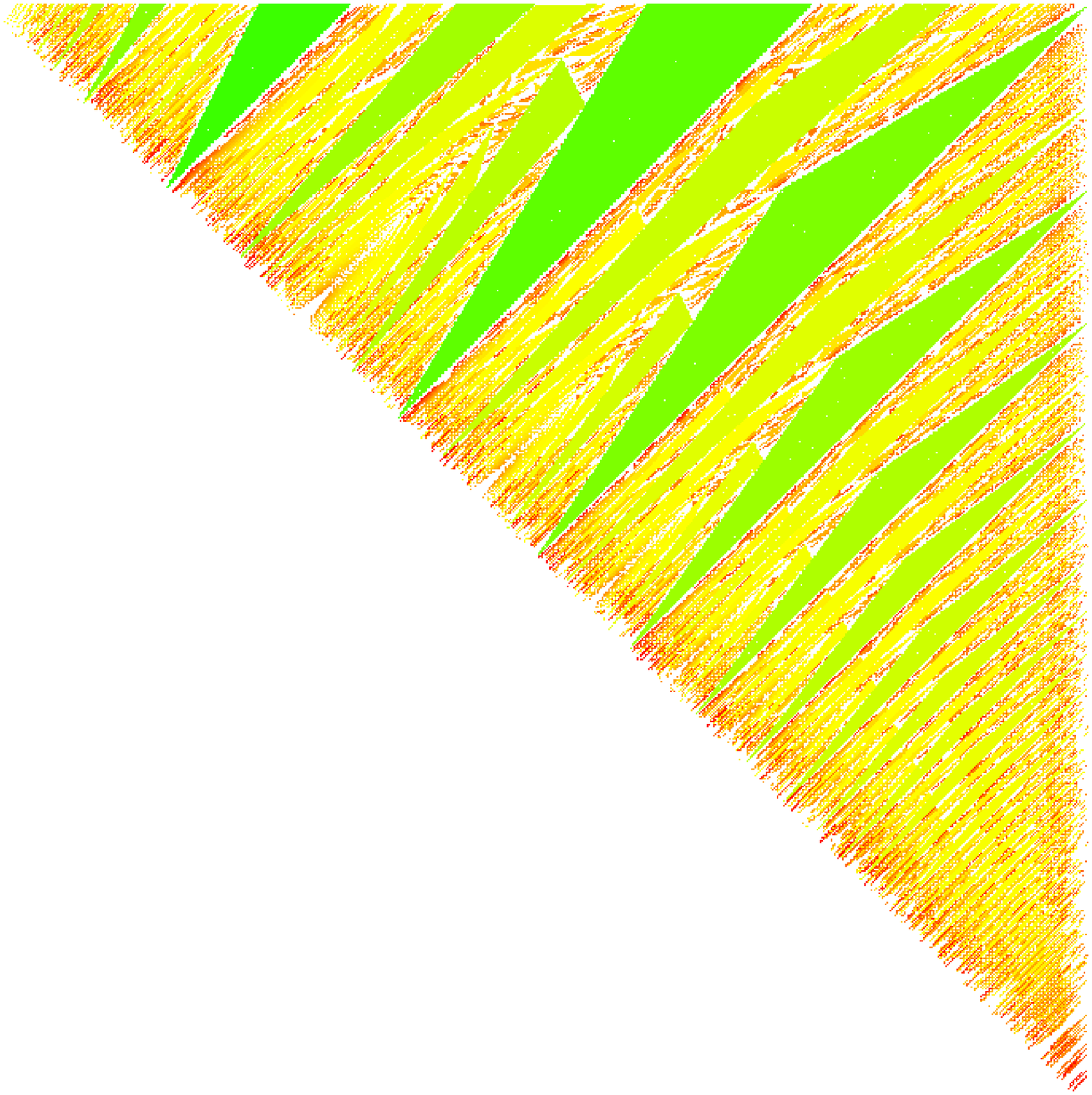}{Detail of $\cS(\cP_0)$ in $[.4,.5]\times[.5,.6]$}
\end{document}